\documentclass[11pt,reqno]{amsart}
\RequirePackage[T1]{fontenc}
\RequirePackage{lmodern}
\RequirePackage{amsmath,amssymb,amsthm,mathtools}
\RequirePackage{microtype}
\RequirePackage{geometry}
\RequirePackage{enumitem}
\RequirePackage{booktabs,array,longtable}
\RequirePackage{xcolor}
\RequirePackage[hidelinks]{hyperref}
\hypersetup{pdfauthor={Jie Wang and Dan Yan},pdftitle={Quartic Keller Maps and Tameness in Three Variables}}
\numberwithin{equation}{section}
\theoremstyle{plain}
\newtheorem{theorem}{Theorem}[section]
\newtheorem{proposition}[theorem]{Proposition}
\newtheorem{lemma}[theorem]{Lemma}
\newtheorem{corollary}[theorem]{Corollary}
\theoremstyle{definition}

\theoremstyle{remark}
\newtheorem{remark}[theorem]{Remark}
\newcommand{\C}{\mathbb C}
\newcommand{\A}{\mathbb A}
\newcommand{\PP}{\mathbb P}
\newcommand{\GA}{\operatorname{GA}}

\newcommand{\Aut}{\operatorname{Aut}}
\newcommand{\tr}{\operatorname{tr}}
\newcommand{\adj}{\operatorname{adj}}
\newcommand{\rank}{\operatorname{rank}}
\newcommand{\Span}{\operatorname{span}}
\newcommand{\Frac}{\operatorname{Frac}}
\newcommand{\Ker}{\operatorname{ker}}

\newcommand{\dd}{\mathrm d}

\newcommand{\J}{\operatorname J}

\newcommand{\nuu}{\boldsymbol\nu}
\newcommand{\XX}{\boldsymbol X}
\newcommand{\Mi}{\mathcal M}
\newcommand{\Tri}{\mathcal T}
\newcommand{\calD}{\mathcal D}
\newcommand{\calE}{\mathcal E}
\newcommand{\br}[2]{\{#1,#2\}}
\allowdisplaybreaks[2]
\setlist{nosep,leftmargin=*}
\title[Quartic Keller Maps and Tameness]{Quartic Keller Maps and Tameness in Three Variables}
\author{Jie Wang}
\address{Jie Wang, State Key Laboratory of Mathematical Sciences, Academy of Mathematics and Systems Science, Chinese Academy of Sciences, Beijing, China}
\email{wangjie212@amss.ac.cn}
\author{Dan Yan}
\address{Dan Yan, MOE-LCSM, School of Mathematics and Statistics,
Hunan Normal University, Changsha, China}
\email{yan-dan-hi@163.com}
\date{\today}
\begin{document}
\begin{abstract}
We study invertibility and tameness for polynomial maps of affine three-space
of total degree at most four. We prove that every complex Keller map in this
degree range is a polynomial automorphism and, independently, that every
complex polynomial automorphism of total degree at most four is tame.
The invertibility argument combines homogeneous factorization of the leading
part, invariant rings of conic pencils, and residue obstructions on cubic
fibers. The tameness argument uses a criterion based on a constant second
directional derivative, low-degree coordinate classifications, and local
tameness of plane polynomial automorphisms over a polynomial parameter ring.
The invertibility conclusion extends to every field of characteristic zero.
\end{abstract}
\maketitle

\section{Introduction}\label{sec:introduction}
Polynomial maps with constant nonzero Jacobian determinant connect local
algebraic properties with global questions about affine space. The chain rule
shows that every polynomial automorphism has constant nonzero Jacobian
determinant. The converse is the assertion traditionally known as the
Jacobian conjecture, originating in Keller's work~\cite{Keller}. Its study
has led to methods involving formal inverses, homogeneous reductions,
polynomial derivations, and the geometry of fibers; see the survey of
Bass, Connell, and Wright~\cite{BCW} and the monographs
\cite{Essen,EKC}.

The general conjecture was disproved in dimension three by the
counterexample announced by Levent Alp\"oge on 19 July 2026~\cite{Alpoge}.
His map has component total degrees $(7,6,4)$ and constant Jacobian
determinant $-2$, yet identifies three distinct points. Its ordinary total
degree is therefore seven, while its geometric degree is three; see
\cite[Section~3.4]{Gao}. Tao~\cite{Tao} gives a geometric explanation of
the example, and Gao~\cite{Gao} develops the tangent-sweep construction
and higher-dimensional extensions. Van den Essen~\cite{E26} gives an
elementary construction recovering Alp\"oge's example up to linear
coordinate changes. These examples demonstrate the role
of behavior at infinity: a constant nonzero Jacobian determinant gives
local invertibility, but does not ensure properness or global injectivity.
They also motivate identifying degree ranges and structural hypotheses
under which polynomial invertibility still follows.

Degree bounds provide a natural setting in which to study this question.
Quadratic Keller maps are invertible in every dimension~\cite{Wang}, and
Vistoli~\cite{Vistoli} gave a conceptual proof of invertibility in dimension
three and degree three. The previously established results discussed here
leave open the invertibility of arbitrary Keller maps in dimension three
of ordinary total degrees four, five, and six.
In dimension two, Moh~\cite{Moh} established
invertibility through degree one hundred; Guccione, Guccione, Horruitiner,
and Valqui~\cite{GGHV} subsequently raised the lower bound for a possible
counterexample to one hundred and eight. The homogeneous cubic reduction
of Bass, Connell, and Wright~\cite{BCW} concerns maps in varying dimension.
It therefore does not replace the analysis of arbitrary maps of a fixed
degree in three variables. Recent work of de Bondt~\cite{dB26} makes
this distinction concrete: by economical stabilization of Alp\"oge's
example, he constructs a cubic Keller counterexample in dimension twelve
and a cubic-linear counterexample in dimension ninety, as well as
counterexamples to forms of Zhao's vanishing conjecture.
These constructions make the increase in dimension explicit when a
high-degree map is replaced by a low-degree one; degree bounds in a fixed
dimension require separate arguments.

Recent algebraic work on the planar conjecture gives further constraints
on possible counterexamples. Lee and Li~\cite{LL24} develop generalized
Magnus expansions and restrict the northeastern vertices of the Newton
polygons of the associated inner polynomials. Valqui and
Ram\'irez~\cite{VR25} compute a Gr\"obner basis and determine the solution
set of a polynomial system arising in a Laurent-series formulation of
the planar problem. These works illustrate how expansions and elimination
constrain Jacobian pairs; our homogeneous coefficient analysis applies
related algebraic tools to quartic leading parts.

The geometry of fibers and their behavior at infinity provide a
complementary perspective. El Hilany's survey~\cite{EH25} places the
Jacobian problem within the topological classification of polynomial
maps, including the study of nonproperness sets and singularities at
infinity. In the real planar setting, Braun, Fernandes,
Gwo\'zdziewicz, and Or\'efice-Okamoto~\cite{BFGO26} prove injectivity
for polynomial maps with nowhere-vanishing Jacobian determinant when
one component has degree at most six, and construct a noninjective
example with a component of degree seven. Their sharp result concerns
real injectivity under a nowhere-vanishing determinant hypothesis.
Our use of cubic components and boundary residues addresses polynomial
invertibility over $\C$ under the constant-Jacobian condition.

Formal inverses also enter recent formulations over coefficient rings.
Hamada, Kato, and Komiya~\cite{HKK25} relate the Jacobian statement for
Tate algebras over a ring with Hausdorff adic topology to the polynomial
statement over its residue ring, and give a non-Archimedean
reformulation over $\mathbb C_p$. Their use of uniqueness of formal
inverses offers a related perspective on coefficient descent, which
appears below in the extension of invertibility to characteristic-zero
fields.

A second question asks how polynomial automorphisms are generated.
An automorphism is called tame if it
is a composition of affine and elementary automorphisms. Jung~\cite{J}
and van der Kulk~\cite{K} proved that every plane polynomial automorphism
is tame. In three variables, Shestakov and Umirbaev~\cite{SU} proved that
Nagata's automorphism is wild, so polynomial invertibility alone does not
settle the tameness question. Blanc and van Santen~\cite{BVS} classified
automorphisms of affine three-space of total degree at most three up to
affine changes of source and target. In particular, every automorphism in
that degree range is tame. Their classification also supplies explicit
normal forms for the low-degree coordinates that occur below.

In related work~\cite{WY26}, we study polynomial maps $H$ in three
variables with nilpotent Jacobian over a field $K$ of characteristic zero.
The classification proposed there reduces normalized maps $H$, with
$H(0)=0$ and linearly independent components, by constant linear
conjugation to the family
\[
\mathcal N_P(x,y,z)
=\bigl(P(y+x^2),\ z-2xP(y+x^2),\ P(y+x^2)^2\bigr),
\qquad P\in K[t]\setminus K,\quad P(0)=0.
\]
This family has total degree $4\deg P$, so its first nontrivial members
are quartic. The proposed classification leads to explicit inverses and
tame factorizations for the one-parameter families $\mathrm{id}+sH$, together with
stable-tameness results for specified higher-dimensional block extensions.
The methods use algebraically dependent constant linear combinations of
the components, generic-fiber geometry, and polynomial derivations.
These tools are also relevant here, but the present invertibility proof
uses the total degree bound directly and is independent of that
classification. In particular, no nilpotence assumption on
$\J(F-\mathrm{id})$ is imposed.

We address invertibility and tameness in total degree at most four. The invertibility
argument applies to arbitrary Keller maps of total degree at most four,
without an initial automorphism assumption. The tameness argument assumes
polynomial invertibility and studies elementary factorizations in the same
degree range. Their conclusions combine as follows.

\begin{theorem}\label{thm:combined}
Let $F=(F_1,F_2,F_3)\in\C[x,y,z]^3$ and suppose
\[
  \max_{1\le i\le3}\deg F_i\le4,
  \qquad \det\J F\in\C^*.
\]
Then $F$ is a tame polynomial automorphism of $\A^3_\C$.
\end{theorem}

The degree hypothesis is essential to the scope of this statement.
Throughout the paper, degree means ordinary total degree. The geometric
degree of a dominant generically finite map instead counts the points in a
generic fiber. Alp\"oge's counterexample has total degree seven, and
Gao's three-dimensional example of geometric degree four has component
total degrees $(4,11,12)$~\cite[Theorem~3.5]{Gao}. Both lie outside the
total degree range of Theorem~\ref{thm:combined}. The theorem identifies
a range in which the Keller condition still forces invertibility and
the resulting automorphisms admit tame factorizations.

The proof of invertibility begins with the homogeneous quartic part of
$F$. Let $\rho$ be the dimension of the span of its three components.
When $\rho=3$, homogeneous factorization reduces the leading part to four
degree patterns. The binary and common-factor patterns are controlled by
successive homogeneous coefficients of the Jacobian determinant. The
remaining pattern is a Veronese map composed with a pencil of quadrics.
The invariant ring of the associated Jacobian derivation restricts this
pencil, and explicit coefficient calculations treat the resulting normal
forms. This part of the argument either excludes a leading form or produces
a constant second directional derivative.

When $\rho\le2$, a linear target change gives a component $R$ of degree
at most three. On every smooth fiber of $R$, the Keller identity makes
the residue two-form exact. Periods along the boundary at infinity restrict
the cubic leading part of $R$, while a second residue calculation excludes
the cuspidal case and controls a conic tangent to a line. The remaining
components admit coordinate completions with controlled inverse degree,
apart from the scalar $x+x^2y$. A formal Laurent-residue argument and a
leading-coefficient obstruction exclude that scalar. The coordinate
completions reduce invertibility to plane Keller maps of degree at most
twenty, well within the established planar bounds~\cite{Moh,GGHV}.

The tameness proof uses the preceding leading-part calculations together
with tools specific to automorphisms. We first establish a criterion,
without a degree bound, based on a constant second directional derivative.
Its affine case uses de Bondt's classification of quasi-translations
\cite{QT}; its quadratic case uses irreducibility of coordinate fibers.
For maps with a component of degree at most three, the classification
of Blanc and van Santen~\cite{BVS} reduces the problem to a small number
of coordinate forms. Components of the form $x+C(y,z)$ lead to families
of plane automorphisms over $\C[t]$. Furter's work on fixed multidegrees
\cite{F} and the local-to-global theorem of Berson, Dubouloz, Furter,
and Maubach~\cite{BDFM} provide the relevant tameness criteria.
The remaining cubic form is treated by homogeneous Jacobian identities.

Section~\ref{sec:preliminaries} fixes the notation and records the
algebraic inputs and homogeneous factorization lemma.
Section~\ref{sec:invertibility} proves invertibility and gives the explicit
coefficient certificates used for the conic pencils.
Section~\ref{sec:tameness} proves tameness, reusing the leading-part
reductions from Section~\ref{sec:invertibility}.
Section~\ref{sec:conclusion} discusses the scope of the results and
directions for further work.

\section{Notation and preliminaries}\label{sec:preliminaries}
We collect the notation and algebraic facts used in both proofs. Except in
Corollary~\ref{cor:charzero}, the base field is $\C$.

\subsection{Polynomial maps and degree conventions}
The degree of a polynomial map is its \emph{ordinary total degree}
$\deg F=\max_i\deg F_i$.
A polynomial map with constant nonzero Jacobian determinant is called a
\emph{Keller map}. A polynomial $r$ is a \emph{coordinate} if it is a
component of a polynomial automorphism; such an automorphism is called a
\emph{coordinate completion} of $r$. An elementary automorphism fixes all
but one variable and adds to that variable a polynomial in the others.
The group generated by affine and elementary automorphisms is the tame
automorphism group. For a coefficient ring $R$, the notation $\GA_2(R)$
denotes the group of $R$-automorphisms of $R[u,v]$.
We regard polynomial maps as column vectors and write
$\J F=(F_x,F_y,F_z)$ for the Jacobian matrix, whose columns are the
partial derivatives of $F$.

Independent affine changes of source and target preserve the Keller
condition, the total degree bound, and tameness. They need not be simultaneous
conjugations. A normalization of the linear part will therefore be used
together with a leading normal form only when their compatibility has been
justified. The term quartic will refer to the degree range at most four.
For a nonzero constant vector $e$, $\partial_e$ denotes differentiation in
the direction $e$. All linear independence statements about polynomial
components are over $\C$ unless another field is specified.

\subsection{Homogeneous Jacobian identities}

After subtracting the constant term, write
\[
F=H_4+H_3+H_2+H_1,\qquad J_i=\J H_i,
\]
where $H_i$ is homogeneous of degree $i$ and a missing homogeneous part is zero. We denote the constant Jacobian determinant by $\kappa\ne0$. The symbol
$\delta$ used later denotes a coefficient parameter.
For three-by-three matrices set
\[
\Mi(A,B)=[t^2]\det(A+tB),\qquad
\Tri(A,B,C)=[stu]\det(sA+tB+uC).
\]
Thus $\Mi$ has one row from $A$ and two rows from $B$, and $\Tri$ is the sum of the six row allocations. If $\Delta_j$ denotes the homogeneous degree-$j$ part of $\det\J F$, then
\begin{align}
\Delta_9&=\det J_4,\notag\\
\Delta_8&=\tr(\adj J_4\,J_3),\notag\\
\Delta_7&=\tr(\adj J_4\,J_2)+\Mi(J_4,J_3),\label{eq:degree-identities}\\
\Delta_6&=\tr(\adj J_4\,J_1)+\Tri(J_4,J_3,J_2)+\det J_3,\notag\\
\Delta_5&=\Mi(J_4,J_2)+\Tri(J_4,J_3,J_1)+\Mi(J_2,J_3).\notag
\end{align}
For a Keller map, all the expressions in \eqref{eq:degree-identities} vanish.

\subsection{Algebraic facts}
The following results will be used throughout the proof.

A Keller map $\A^2_\C\to\A^2_\C$ of degree at most twenty is a polynomial automorphism. This is a much weaker consequence of Moh's bound of one hundred~\cite{Moh}.

An injective Keller map of complex affine space is a polynomial
automorphism~\cite{Ax,BCW}. Indeed, injectivity implies surjectivity by
Ax's theorem, and a bijective \'etale morphism in this setting is an
isomorphism.

\subsection{Closed polynomial generators}\label{sec:closed-generators}
We give proofs of the facts about closed polynomials and algebraic
dependence needed below. These are standard facts; see \cite{AP} for
background and further characterizations. All subsequent uses of these
facts refer to the proofs below.

Write $A=\C[x_1,\ldots,x_n]$ and $L=\Frac A$. A nonconstant
polynomial $s\in A$ is called \emph{closed} if every element of $A$
integral over $\C[s]$ belongs to $\C[s]$. A subfield $E\subset L$
is \emph{relatively algebraically closed} in $L$ if every element of
$L$ algebraic over $E$ belongs to $E$.

For any nonconstant polynomial $s$ one has
\begin{equation}\label{eq:intersection}
\C(s)\cap\C[x_1,\ldots,x_n]=\C[s].
\end{equation}
Indeed, a B\'ezout identity for coprime univariate polynomials $a,b$ remains a B\'ezout identity for $a(s),b(s)$. If $a(s)/b(s)$ is polynomial, $b(s)$ is a unit, and therefore $b$ is constant. The same argument shows that if $a(s)/b(s)=f/g$, where $a,b$ are coprime
and $f,g$ are polynomials, then $b(s)$ divides $g$.

\begin{lemma}[Normal subalgebras in one variable]\label{lem:normal-subalgebra}
Let $B$ be a $\C$-subalgebra of $\C[t]$ that contains a nonconstant
polynomial and is integrally closed in its fraction field.
Then $B=\C[u]$ for some $u\in\C[t]$.
\end{lemma}
\begin{proof}
Put $E=\Frac B$, and choose a monic $u\in B$ of least positive degree
$m$. Since $t$ satisfies the monic equation $u(T)-u(t)=0$ over $B$,
$\C[t]$ is integral over $B$. Consequently
$E\cap\C[t]=B$.

Let $M(T)\in E[T]$ be the monic minimal polynomial of $t$ over $E$.
It divides $u(T)-u(t)$ in $\C(t)[T]$. Both $M$ and the monic quotient
$N$ lie in $\C[t,T]$: their coefficients are integral over $\C[t]$
and belong to $\C(t)$, and $\C[t]$ is integrally closed. In particular,
every coefficient of $M$ lies in $E\cap\C[t]=B$. Some coefficient
of $M$ is nonconstant in $t$, since $M(t)=0$ and $t$ is transcendental
over $\C$. Its degree is at least $m$, by the choice of $u$. Taking
total degrees in the two independent variables $t,T$ gives
\[
 m=\deg_{t,T}\bigl(u(T)-u(t)\bigr)
   =\deg_{t,T}M+\deg_{t,T}N.
\]
Thus $\deg_{t,T}M\ge m$ forces $N$ to be constant, and
$[\C(t):E]=\deg_T M=m$. The equation $u(T)-u(t)=0$ also gives
$[\C(t):\C(u)]\le m$. Since $\C(u)\subset E$, the tower law
yields $E=\C(u)$. Finally, the one-variable instance
of \eqref{eq:intersection} gives $B=E\cap\C[t]=\C[u]$.
\end{proof}

\begin{proposition}[Closed generators]\label{prop:closed-generators}
The following statements hold in $A$.
\begin{enumerate}[label=\textup{(\roman*)}]
\item Every nonconstant $r\in A$ has a representation $r=h(s)$,
where $h\in\C[T]$ and $s\in A$ is closed.
\item A nonconstant $s\in A$ is closed if and only if $\C(s)$ is
relatively algebraically closed in $L$.
\item If $p\in A$ is nonconstant and $q\in A$ is algebraic over
$\C(p)$, then $p=P(s)$ and $q=Q(s)$ for a common closed polynomial
$s\in A$ and some $P,Q\in\C[T]$.
\end{enumerate}
\end{proposition}
\begin{proof}
For (i), let $B$ be the integral closure of $\C[r]$ in $L$.
Every element of $B$ is integral over $A$ and hence belongs to $A$,
since $A$ is integrally closed. Transitivity of integrality also
shows that $B$ is integrally closed in its fraction field.

Choose an affine line $\ell:\A^1\to\A^n$ on which $r$ is
nonconstant. Such a line can be chosen through two points with
distinct $r$-values. Restriction to $\ell$ defines a homomorphism
$B\to\C[t]$, and we claim that it is injective. Indeed, any
nonzero $b\in B$ satisfies a monic equation over $\C[r]$ whose
constant term $c(r)$ is nonzero: remove any common power of the
indeterminate from an integral equation for $b$. This equation gives
$c(r)\in bB$. If $b$ vanished identically on $\ell$, then $c(r)$ would also vanish
there, contradicting the nonconstancy of $r|_\ell$.
Lemma~\ref{lem:normal-subalgebra} applied to the image therefore
gives $B=\C[s]$ for some $s\in A$. In particular $r=h(s)$.
If an element of $A$ is integral over $\C[s]=B$, it is integral
over $\C[r]$ and thus belongs to $B$. Hence $s$ is closed.

For (ii), suppose first that $s$ is closed and that $v\in L$ is
algebraic over $\C(s)$. Clearing denominators gives an equation
\[
 c_d(s)v^d+c_{d-1}(s)v^{d-1}+\cdots+c_0(s)=0,
 \qquad c_i\in\C[T],\quad c_d\ne0.
\]
Multiplication by $c_d(s)^{d-1}$ shows that $c_d(s)v$ is integral
over $\C[s]$. It belongs to $A$ because $A$ is integrally closed,
and then to $\C[s]$ because $s$ is closed. Thus $v\in\C(s)$.
Conversely, if $\C(s)$ is relatively algebraically closed in $L$,
every element of $A$ integral over $\C[s]$ belongs to $\C(s)$,
and hence to $\C[s]$ by \eqref{eq:intersection}.

For (iii), apply (i) to $p$ to obtain $p=P(s)$ with $s$ closed.
The polynomial $q$ is algebraic over $\C(s)$, so (ii) gives
$q\in\C(s)$. Equation~\eqref{eq:intersection} then gives
$q\in\C[s]$.
\end{proof}

\begin{lemma}[Differential criterion and constants]\label{lem:differential-constants}
Let $r\in L\setminus\C$ and $v\in L$. Then $v$ is algebraic
over $\C(r)$ if and only if
\[
 (\partial_{x_i}r)(\partial_{x_j}v)
 -(\partial_{x_j}r)(\partial_{x_i}v)=0
 \qquad(1\le i<j\le n).
\]
In particular, if $s\in\C[x,y]$ is closed and
$D_s(v)=\{s,v\}=s_xv_y-s_yv_x$, then
\[
 \ker\bigl(D_s:\C(x,y)\to\C(x,y)\bigr)=\C(s),
 \qquad
 \ker\bigl(D_s:\C[x,y]\to\C[x,y]\bigr)=\C[s].
\]
\end{lemma}
\begin{proof}
If $v$ is algebraic over $\C(r)$, differentiate its minimal
polynomial relation. Characteristic zero makes the derivative of that
polynomial with respect to its indeterminate nonzero at $v$. Thus
$\partial_{x_i}v=a\,\partial_{x_i}r$ for all $i$, with the same
$a\in L$, and the displayed minors vanish.

Conversely, if $r,v$ are algebraically independent, extend them to a
transcendence basis of $L/\C$. The partial derivations with respect
to $r$ and $v$ on the resulting purely transcendental subfield extend
to $\C$-derivations $\delta,\epsilon$ of $L$, since the remaining
extension is finite and separable. Explicitly, for an algebraic
element $\alpha$ with minimal polynomial $P$, the extension is
determined by $\delta(\alpha)=-\delta(P)(\alpha)/P'(\alpha)$,
where $\delta(P)$ differentiates the coefficients; the same formula
applies to $\epsilon$. Hence
\[
 \delta(r)=1,\quad\delta(v)=0,\qquad
 \epsilon(r)=0,\quad\epsilon(v)=1.
\]
Every $\C$-derivation $\eta$ of $L$ has the form
$\sum_i\eta(x_i)\partial_{x_i}$. Vanishing of all the displayed
minors would therefore imply
$\delta(r)\epsilon(v)-\epsilon(r)\delta(v)=0$, a contradiction.
The two kernel identities now follow from
Proposition~\ref{prop:closed-generators}(ii) and
\eqref{eq:intersection}.
\end{proof}

\subsection{Homogeneous factorization}
The next lemma organizes the possible quartic leading parts when their
components are linearly independent.
\begin{lemma}\label{lem:factor}
Let $H=(H_1,H_2,H_3)$ consist of linearly independent homogeneous forms of a common positive degree $d$, and suppose $\det\J H=0$. There are coprime homogeneous forms $p,q$ of the same degree $m\ge1$, coprime binary forms $\Phi_i$ of a common degree $e\ge2$, and a homogeneous polynomial $g$ of degree $r\ge0$ such that
\begin{equation}\label{eq:factor}
H_i=g\Phi_i(p,q),\qquad d=r+em.
\end{equation}
For $d=4$ the possibilities for $(r,m,e)$ are
\begin{equation}\label{eq:patterns}
(0,1,4),\quad(1,1,3),\quad(2,1,2),\quad(0,2,2).
\end{equation}
\end{lemma}
\begin{proof}
Euler's identity shows that the radial value vector belongs to the column space of $\J H$. The projective differential rank is one less than the affine rank. If the latter were at most one, the projective image would be a point, contrary to independence. Thus $\rank\J H=2$, and the image of $[H_1:H_2:H_3]$ is a curve.
The restriction to a suitable line in the source is nonconstant. It gives a dominant rational map from $\PP^1$ to the normalization of that curve, which is therefore rational in characteristic zero. Choose a rational parameter on the normalization; its pullback is $p/q$ for coprime homogeneous forms of equal degree $m$. The normalization map to $\PP^2$ is represented by binary forms $\Phi_i$ of a common degree $e$ without common projective zero. Since the image is not a line, $e\ge2$.

The polynomials $\Phi_i(p,q)$ are coprime. Indeed, powers of both binary variables belong to the ideal generated by the $\Phi_i$, so a common irreducible divisor of the substituted forms would divide both $p$ and $q$. The rational proportionality factor between $H_i$ and $\Phi_i(p,q)$
is consequently a homogeneous polynomial $g$. This proves \eqref{eq:factor}; elementary degree comparison gives \eqref{eq:patterns}.
\end{proof}

\section{Invertibility of quartic Keller maps}\label{sec:invertibility}
This section proves the invertibility statement underlying
Theorem~\ref{thm:combined}. The proof separates maps with three independent
quartic leading components from those with a component of degree at most
three. The former are treated by homogeneous Jacobian identities; the
latter by residue obstructions and coordinate completions with controlled
inverse degree.

\begin{theorem}\label{thm:main}
Let $F\in\C[x,y,z]^3$ satisfy
\[
\deg F:=\max_i\deg F_i\le4,
\qquad \det\J F\in\C^*.
\]
Then $F$ is a polynomial automorphism of $\A^3_\C$.
\end{theorem}

\subsection{A constant second directional derivative}\label{sec:direction}
We begin with an invertibility criterion that will apply to several leading
forms. The affine reduction below is also used to exclude the exceptional
scalar component in Section~\ref{sec:exceptional-scalar}.

Write $\{f,g\}=f_xg_y-f_yg_x$.
\begin{lemma}[Affine reduction]\label{lem:affine-reduction}
Let $a,b,p,q,r\in\C[x,y]$ and suppose
\[
\det\J(a+zp,b+zq,r)=\kappa\in\C^*.
\]
Then $p,q\in\C[r]$, and $S=pb-qa$ satisfies $\{S,r\}=\kappa$.
\end{lemma}
\begin{proof}
The third column of the Jacobian matrix shows that $p,q$ have no common zero.
By Proposition~\ref{prop:closed-generators}(i), write $r=h(s)$ with
$s$ a closed generator.
Since $h'(s)$ divides every entry of the third Jacobian row, it divides $\kappa$;
thus $h$ is affine linear and $r$ is itself closed.
Interchange the first two components if needed so that $q\ne0$.
The coefficient of $z$ is $p\{q,r\}-q\{p,r\}=0$, so $\{p/q,r\}=0$.
Since $r$ is closed, the rational-function kernel identity in
Lemma~\ref{lem:differential-constants} gives $p/q\in\C(r)$.
Write $p/q=A(r)/B(r)$ in reduced form. B\'ezout gives $p=\eta A(r)$, $q=\eta B(r)$ for a polynomial $\eta$; the absence of a common zero makes $\eta$ constant.
Thus $p,q\in\C[r]$, and the constant coefficient gives
\[
\{pb-qa,r\}=p\{b,r\}-q\{a,r\}=\kappa.
\]
\end{proof}

\begin{proposition}\label{prop:direction}
Let $F$ be a Keller map of degree at most four. If $\partial_e^2F\in\C^3$ for some nonzero constant vector $e$, then $F$ is a polynomial automorphism.
\end{proposition}
\begin{proof}
After a linear source change,
\begin{equation}\label{eq:direction-form}
F=A(x,y)+zB(x,y)+z^2E,\qquad
\deg A\le4,\quad\deg B\le3,\quad E\in\C^3.
\end{equation}
Both cases use the following elementary observation over $k=\C(x)$:
\begin{equation}\label{eq:product-obstruction}
u\in k[y],\quad H\in k[T],\quad\deg H\ge1
\quad\Longrightarrow\quad u_yH(u)\notin k^*.
\end{equation}
Indeed, a nonzero constant product makes $H(u)$ a unit, forcing $u\in k$ and hence $u_y=0$.

\subsubsection*{The affine case \texorpdfstring{$E=0$}{E=0}}
The components of $B=F_z$ have no common zero. The coefficient of $z^2$ in $\det\J F$ is
\begin{equation}\label{eq:affine-B}
\det(B_x,B_y,B)=0.
\end{equation}
Suppose first that the components of $B$ are independent. Homogenize them to their maximum degree $d\le3$. Euler's identity and \eqref{eq:affine-B} show that the homogeneous Jacobian is singular. Lemma~\ref{lem:factor} applies, and $e\ge2$, $d\le3$ force $m=1$. After dehomogenization, the common factor is constant because $B$ has no common zero. The two resulting affine-linear forms have no common zero either, so their gradients are dependent. After an affine change in $x,y$, every component of $B$ depends on $x$ alone.

The coefficient of $z$ in $\det\J F$ now gives
\[
\det(B',A_y,B)=0.
\]
Over $k=\C(x)$, the independent components of $B$ have nonzero Wronskian
\[
\Omega=\det(B'',B',B).
\]
To verify this, use constant row operations to give the three polynomials distinct degrees: the leading Wronskian coefficient is the product of their leading coefficients and a nonzero Vandermonde determinant in those degrees. In particular $B,B'$ are $k$-independent. Expressing $A_y$ in their $k$-span and integrating with respect to $y$ gives
\[
A=uB+vB'+C(x),\qquad u,v\in k[y],\quad C\in k^3.
\]
Direct expansion of $\det(A_x+zB',A_y,B)$ gives
\begin{equation}\label{eq:affine-wronskian}
\det\J F=v_y\bigl(\Omega v+\Lambda\bigr),
\qquad\Lambda=\det(C',B',B).
\end{equation}
Since $\Omega\ne0$, \eqref{eq:product-obstruction} contradicts $\det\J F=\kappa$. Thus the components of $B$ are dependent.

Make a constant target change to write
\begin{equation}\label{eq:affine-dependent}
F=(a+zp,b+zq,r),
\end{equation}
where $a,b,r\in\C[x,y]$ have degree at most four and $p,q$ have degree at most three.
Lemma~\ref{lem:affine-reduction} gives $p,q\in\C[r]$ and a planar Keller map $(S,r)$ with $S=pb-qa$.
Its degree is at most seven, so the planar invertibility theorem applies.
Two points with the same $F$-image have the same values of $r$ and $S=pF_2-qF_1$, hence the same $(x,y)$. Since $p,q$ have no common zero, the common $F$-image then determines $z$.
Thus $F$ is injective and hence an automorphism.

\subsubsection*{The quadratic case \texorpdfstring{$E\ne0$}{E not equal to 0}}
A target change puts \eqref{eq:direction-form} in the form
\begin{equation}\label{eq:quadratic-original}
F=(a+zp,b+zq,c+zr+z^2).
\end{equation}
Here $\deg a,\deg b,\deg c\le4$ and $\deg p,\deg q,\deg r\le3$. The polynomials $p,q$ have no common zero: at a common zero choose $z=-r/2$ to make $F_z=0$.
For coefficient calculations only, complete the square in $z$, obtaining $(a+zp,b+zq,c+z^2)$. This changes $a,b,c$ and leaves $p,q$ fixed; the argument uses no degree
estimate for the modified $a,b,c$. The determinant coefficients are
\begin{align}
[z^3]\det\J F&=2\{p,q\},\notag\\
[z^2]\det\J F&=2(\{p,b\}+\{a,q\}),\notag\\
[z]\det\J F&=2\{a,b\}+p\{q,c\}-q\{p,c\},\label{eq:quadratic-coeff}\\
[1]\det\J F&=p\{b,c\}-q\{a,c\}.\notag
\end{align}
We claim that a constant linear combination of $p,q$ is a nonzero constant.
Suppose not. Then $1,p,q$ are linearly independent: a proportional pair
without common zeros would be constant. The first equation in
\eqref{eq:quadratic-coeff} gives $\{p,q\}=0$. By
Lemma~\ref{lem:differential-constants} and
Proposition~\ref{prop:closed-generators}(iii), there is a common closed
polynomial generator $s$
with $p=P(s)$, $q=Q(s)$. At least one outer polynomial has degree at least
two. The degree bound on the original $p,q$ therefore forces $\deg s=1$.
Choose $s=x$. Constant linear combinations give
\[
1\le m=\deg p<n=\deg q\le3.
\]
Set
\[
W=pq'-qp',\qquad K=p''q'-p'q''.
\]
Both are nonzero, with respective degrees $m+n-1$ and $m+n-3$. The $z^2$ equation gives $p'b_y-q'a_y=0$, so over $\C(x)$
\[
a=p'u,\qquad b=q'u+\beta(x),\qquad u\in\C(x)[y].
\]
The $z$ equation is
\[
Wc_y=-2Kuu_y+2p'\beta'u_y.
\]
Integrating in $y$ gives
\[
c=-\frac KWu^2+h_1u+h_0,\qquad h_0,h_1\in k.
\]
Since $pb_x-qa_x=W'u+Wu_x+p\beta'$ and $pb_y-qa_y=Wu_y$, the $u_x$ terms cancel in the constant determinant equation, giving
\begin{equation}\label{eq:quadratic-wronskian}
\kappa=u_yH(u),\qquad
H(T)=\left(K'-3K\frac{W'}W\right)T^2+\ell_1T+\ell_0,
\end{equation}
where $\ell_0,\ell_1\in k$. The quadratic coefficient is nonzero: otherwise $(K/W^3)'=0$, contrary to the degrees of $K,W$.
Thus \eqref{eq:product-obstruction} again contradicts the Keller condition, proving the claim.

Return to the coordinates in \eqref{eq:quadratic-original}, before the
completion of the square. A constant target change and rescaling make the first component $z+a_0(x,y)$ with $\deg a_0\le4$. Use $t=z+a_0$ as a polynomial source coordinate. For every fixed $t$, the other two components form a planar Keller map with component degrees at most seven and eight: substitution only introduces products of $a_0$ with degree-three coefficients and $a_0^2$. The planar invertibility theorem makes each fiber map injective. Thus $F$ is injective, and hence an automorphism.
\end{proof}

\subsection{Independent leading parts: binary forms and common factors}\label{sec:leading-simple}
We next analyze the patterns in Lemma~\ref{lem:factor} for which the
leading part factors through binary forms. The first lemma controls the
normal vector to their image; the next two translate this control into
homogeneous Jacobian obstructions.
\begin{lemma}\label{lem:normal-annihilator}
If $P=(P_1,P_2,P_3)$ consists of independent binary homogeneous forms of a common positive degree, then $n=P_x\times P_y$ has no nonzero constant annihilator: $n\cdot c=0$, $c\in\C^3$, implies $c=0$.
\end{lemma}
\begin{proof}
A target change sending a nonzero $c$ to $e_3$ would give zero planar Jacobian for the first two components. Homogeneous binary forms of equal degree with zero Jacobian are proportional: on $x\ne0$ write them as $x^df(y/x),x^dg(y/x)$ and obtain $fg'-f'g=0$. This contradicts independence.
\end{proof}

\begin{lemma}\label{lem:binary-leading}\label{lem:affine-leading}
Suppose that, for $s\in\{0,1\}$,
\[
H_4=z^sP_{4-s}(x,y)+R_4(x,y),
\]
where $P_{4-s}$ has linearly independent components and $R_4=0$ when $s=0$.
Then $F_{zz}$ is constant and $F$ is an automorphism.
\end{lemma}
\begin{proof}
Write
\[
H_3=B_3(x,y)+zB_2(x,y)+z^2B_1(x,y)+z^3c,
\qquad n=(P_{4-s})_x\times(P_{4-s})_y.
\]
The highest relevant coefficients of $\Delta_8$ are
\[
[z^{2s+2}]\Delta_8=3n\cdot c,
\qquad
[z^{2s+1}]\Delta_8\big|_{c=0}
=\frac{8-3s}{4-s}\,n\cdot B_1.
\]
For $s=0$ these follow directly from the third row of the adjugate.
For $s=1$, Euler's identity gives the leading terms of its rows as
$-xz n^T/3$, $-yz n^T/3$, and $z^2n^T$; the binary remainder contributes only lower powers.
Lemma~\ref{lem:normal-annihilator} gives $c=0$, and also $B_1=0$ if $B_1$ has rank at most one.
If $B_1$ has rank two, a constant target change makes it $(x,y,0)$.
The coefficient of $z^{s+3}$ in $\Delta_7$ is then
\begin{align*}
&2\det((P_{4-s})_x,e_2,(x,y,0))
+2\det(e_1,(P_{4-s})_y,(x,y,0))
+s\det(e_1,e_2,P_{4-s})\\
&\hspace{25mm}=-(8-3s)(P_{4-s})_3,
\end{align*}
again by Euler's identity. Its vanishing contradicts independence.
Thus $B_1=0$, and $F_{zz}=(H_2)_{zz}$ is constant.
Apply Proposition~\ref{prop:direction}.
\end{proof}

\begin{lemma}\label{lem:common-factor}
No quartic Keller map has leading part
\[
H_4=g(x,y,z)(x^2,xy,y^2),\qquad \deg_zg=2,
\]
where $g$ is homogeneous quadratic.
\end{lemma}
\begin{proof}
Normalize the $z^2$ coefficient of $g$ to one by rescaling the target. Write
\[
N=(x^2,xy,y^2)^T,\quad n=(y^2,-2xy,x^2)^T,
\quad k=g_z\XX-4ge_3.
\]
Direct multiplication gives
\begin{equation}\label{eq:common-adj}
\adj\J(gN)=-gkn^T.
\end{equation}
For $S=n^TH_3$, Euler's identity turns $\Delta_8=0$ into
\[
4gS_z=3g_zS.
\]
If $S\ne0$, comparison over $\C(x,y)[z]$ gives $4\deg_zS=6$, which is impossible. Thus $S=0$, whose elementary syzygies give
\[
H_3=(2xU,yU+xV,2yV)
\]
with $U,V$ homogeneous quadratic. A source translation in $x,y$ removes their $z^2$ coefficients: the $z^2$ coefficient of $\partial_{(a,b,0)}H_4$ is $(2ax,ay+bx,2by)$. The new $H_3$ is affine in $z$.
If $c$ is the constant $z^2$ coefficient of $H_2$, the coefficient of $z^5$ in $\Delta_7$, from \eqref{eq:common-adj}, is $4n^Tc$. Thus $c=0$ and $H_2$ is affine in $z$. The coefficient of $z^4$ in $\Delta_6$ is now $2n^T(H_1)_z$; all the other degree-six terms have $z$-degree at most three. Therefore $(H_1)_z=0$.
The derivative $F_z$ consequently vanishes identically on $x=y=0$: its contributions from $H_4,H_3,H_2$ contain respectively $N$, a binary quadratic, and a binary linear form. This contradicts the Keller condition.
\end{proof}

\subsection{Coprime conic pencils without a square}\label{sec:nonsquare}
The remaining pattern in \eqref{eq:patterns} factors through a pencil of
quadrics. Write $\nuu(s,t)=(s^2,st,t^2)^T$. We first determine the polynomial
first integrals of the pencil derivation, using homogeneous degrees on its
generic conic.

\begin{lemma}[Smooth generic member]\label{lem:smooth-generic-conic}
Let $p,q\in\C[x,y,z]$ be coprime homogeneous quadrics whose pencil contains no square of a nonzero linear form. Then $p,q$ are algebraically independent, and the generic projective member of their pencil is a smooth conic.
\end{lemma}
\begin{proof}
The quadrics are linearly independent by coprimeness. Two homogeneous forms of the same degree that are algebraically dependent are proportional: decompose a polynomial relation by degree, and use the factorization of the resulting homogeneous binary relation over $\C$. Thus $p,q$ are algebraically independent.

Suppose every member of the pencil were singular. No nonzero member has rank one, since over $\C$ such a quadratic form is a square up to a nonzero scalar. Choose a rank-two member as $q=xy$ and write
\[
p=ax^2+bxy+cy^2+dxz+eyz+fz^2.
\]
The determinant of the symmetric matrix of $p+txy$ is
\[
acf+\frac{(b+t)de}{4}-\frac{ae^2+cd^2}{4}
       -\frac{f(b+t)^2}{4}.
\]
Its identically vanishing coefficients give $f=0$, $de=0$, and $ae^2+cd^2=0$. If $d\ne0$, then $e=c=0$, so $x$ divides both $p,q$. If $e\ne0$, then $d=a=0$, so $y$ divides both. If $d=e=0$, both quadrics are binary. The homogeneous quadratic determinant of a binary pencil either vanishes identically or has a projective zero over $\C$. In either event, some nonzero binary member has rank one, hence is a square. All cases contradict the hypotheses. The ternary pencil determinant is therefore nonzero, which proves the assertion also at its generic parameter.
\end{proof}

\begin{lemma}\label{lem:pencil-invariants}
Let $p,q\in\C[x,y,z]$ be coprime homogeneous quadrics. If their pencil contains no square of a nonzero linear form, then
\[
\Ker_{\C[x,y,z]}\calD=\C[p,q],\qquad
\calD=(\nabla p\times\nabla q)\cdot\nabla.
\]
\end{lemma}
\begin{proof}
The inclusion $\C[p,q]\subseteq\Ker\calD$ follows from the definition of $\calD$.
Since $\calD$ raises degree by one, it suffices to treat a nonzero homogeneous invariant $f$ of degree $d>0$.

Put $K=\C(\lambda)$, with $\lambda$ indeterminate. Lemma~\ref{lem:smooth-generic-conic} shows that
\[
C=\{p-\lambda q=0\}\subset\PP^2_K
\]
is a smooth conic.
The two projective conics $p=0$ and $q=0$ have a common point over $\C$ by B\'ezout. This is a $K$-rational point of $C$, and projection from it gives $K(C)\simeq K(t)$.
Thus a rational function on $C$ with zero differential belongs to $K$.

Restriction to the affine cone over $C$ does not annihilate a nonzero polynomial with complex coefficients: clear denominators in $\lambda$ and substitute $\lambda=p/q$.
At its generic point, $\nabla p,\nabla q$ are independent, and $\calD f=0$ gives
$\nabla f=u\nabla p+v\nabla q$.
Euler's identity yields $2q(\lambda u+v)=d\,f$.
On a projective chart of $C$ we have $\dd p=\lambda\dd q$, hence
\begin{equation}\label{eq:homogeneous-invariant}
q\,\dd f=\frac d2 f\,\dd q,
\qquad
\dd\left(\frac{f^2}{q^d}\right)=0.
\end{equation}
The degree-zero ratio $f^2/q^d$ is therefore in $K$.
Substituting $\lambda=p/q$ gives
\[
f^2=q^d c(p/q),\qquad c\in\C(\lambda).
\]
The denominator factors $q$ and $p-aq$ for distinct $a\in\C$ are pairwise coprime. Polynomiality of $f^2$ consequently forbids finite poles of $c$ and forces $\deg c\le d$.
Thus
\[
f^2=H(p,q),\qquad H\in\C[s,t]\text{ homogeneous of degree }d.
\]
Distinct binary linear factors of $H$ give coprime pencil members.
Each member is squarefree, since a non-squarefree quadratic over $\C$ is a scalar multiple of a square.
Unique factorization in $\C[x,y,z]$ therefore forces every linear factor of $H$ to have even multiplicity.
Hence $H=G^2$ for a binary form $G$, and $f=\pm G(p,q)\in\C[p,q]$.
This proves the assertion in every degree, including the impossibility of odd-degree invariants.
\end{proof}

\begin{proposition}\label{prop:square-forced}
If a quartic Keller map has $H_4=\nuu(p,q)$ for coprime homogeneous quadrics $p,q$, then the pencil $\langle p,q\rangle$ contains a square.
\end{proposition}
\begin{proof}
Suppose not, and put $n=(q^2,-2pq,p^2)^T$. The universal adjugate identity is
\begin{equation}\label{eq:veronese-adj}
\adj\J\nuu(p,q)=2(\nabla p\times\nabla q)n^T.
\end{equation}
Since $\calD n=0$, it implies
\[
\tr(\adj\J\nuu(p,q)\,\J U)=2\calD(n^TU).
\]
The degree-eight equation gives $\calD(n^TH_3)=0$. The invariant $n^TH_3$ has degree seven, whereas all homogeneous degrees in $\C[p,q]$ are even. Lemma~\ref{lem:pencil-invariants} therefore gives $n^TH_3=0$. Reducing this relation modulo $p$ and $q$ gives
\[
H_3=(2pa,qa+pb,2qb)
\]
for homogeneous linear forms $a,b$. Set $P=p+a$, $Q=q+b$. Absorbing the quadratic terms introduced by completing these products, write
\[
F=\nuu(P,Q)+R_2+L,
\]
where $R_2$ is homogeneous quadratic and $L$ homogeneous linear. The map $\nuu(P,Q)$ has identically singular Jacobian. Consequently its degree-seven determinant equation gives $\calD(n^TR_2)=0$, and hence
\[
n^TR_2=A p^3+B p^2q+C pq^2+D_0q^3
\]
for constants $A,B,C,D_0$.
Reduction modulo $p$ gives $R_{2,1}=D_0q+up$, and reduction modulo $q$ gives $R_{2,3}=Ap+vq$, with constants $u,v$. Substitution then gives
\[
R_{2,2}=\frac{v-B}{2}p+\frac{u-C}{2}q.
\]
Thus $R_2=M(p,q)$ for a constant three-by-two matrix $M$, and
\[
F=\Phi(P,Q)+K\XX,\qquad
\Phi(s,t)=\nuu(s,t)+M(s,t),
\]
where $K$ is constant. Since $\Phi(P,Q)$ still factors through two polynomials, the degree-six equation is now $\calD(n^TK\XX)=0$. This invariant has degree five, so it vanishes. Coprimeness forces $p$ to divide the first linear component of $K\XX$ and $q$ the third. Both components are zero, and the middle is then zero as well. Thus $K=0$, contradicting the nonzero Jacobian of $F$.
\end{proof}

\subsection{The two nonbinary square pencils}\label{sec:square}
Proposition~\ref{prop:square-forced} reduces the conic-pencil analysis to
pencils containing a square. We classify these pencils and then use the
degree-seven and degree-six Jacobian equations to exclude the nonbinary
normal forms.
If the pencil contains a square, choose it as $p=x^2$. Classify the rank of $q$ on the plane $x=0$. Rank two gives $q=yz$ after a linear change of $y,z$, further linear changes to remove the cross terms involving $x$, and subtraction of a multiple of $x^2$. Rank one gives $q=y^2$ or $q=y^2+xz$. Rank zero would make $q$ divisible by $x$, contradicting coprimeness. Pencil changes induce invertible constant changes of the three Veronese components. Hence the possibilities are
\begin{equation}\label{eq:square-normalforms}
(p,q)=(x^2,y^2),\qquad(x^2,yz),\qquad(x^2,y^2+xz).
\end{equation}
The first is covered by Lemma~\ref{lem:binary-leading}.

\subsubsection{Invariant rings and the complete cubic kernel}
For either remaining $q$, let $n=(q^2,-2x^2q,x^4)^T$. The polynomial invariant ring of the Jacobian derivation of $x^2,q$ is $\C[x,q]$. For $q=yz$ this is the monomial weight calculation for $z\partial_z-y\partial_y$. For $q=y^2+xz$, localize at $x$ and use coordinates $x,q,y/x$. The derivation becomes a nonzero multiple of differentiation in the last coordinate, so its kernel in that localized ring is $\C[x,x^{-1},q]$. Intersecting with $\C[x,y,z]$ removes negative powers of $x$: if $x^{-m}A(x,q)$ is polynomial and $A$ is not divisible by $x$, reduction at $x=0$ would give $A(0,y^2)=0$, a contradiction. Thus the intersection is $\C[x,q]$.

The degree-eight equation says that $n^TH_3$ is a degree-seven element of this ring, and therefore a combination of $x^7,x^5q,x^3q^2,xq^3$. Reduction modulo $x^2$ and modulo $q$ gives the complete cubic family
\begin{equation}\label{eq:cubic-kernel}
H_3=(2x^2u+\alpha xq,\ qu+x^2v+\beta_0x^3+\gamma xq,\ 2qv+\delta x^3),
\end{equation}
where $u,v$ are arbitrary homogeneous linear forms. More explicitly, the first reduction gives $H_{3,1}-\alpha xq\in(x^2)$, the second gives $H_{3,3}-\delta x^3\in(q)$, and substitution determines the middle component. The ten parameters in this representation are independent: setting the three components to zero and using $\gcd(x,q)=1$ successively forces $\alpha,u,\delta,v,\beta_0,\gamma$ to vanish.

The quadratic normal forms below come from linear elimination. Define
\[
\calE=z\partial_z-y\partial_y\quad(q=yz),\qquad
\calE=-x\partial_y+2y\partial_z\quad(q=y^2+xz).
\]
Since $\calD=2x\calE$, the coefficient operator on any homogeneous vector $U$ is
\begin{equation}\label{eq:coefficient-operator}
\mathcal L(U)=\tr(\adj J_4\,\J U)=4x\,n^T\calE U.
\end{equation}
On quadratic vectors its kernel is precisely
\begin{equation}\label{eq:quadratic-kernel}
\{M(x^2,q):M\in\operatorname{Mat}_{3\times2}(\C)\}.
\end{equation}
Indeed, $n^TU$ must belong to $\C[x,q]$. In degree six, reduce its four possible monomials modulo $x^2,q$ to obtain \eqref{eq:quadratic-kernel}. Thus the quadratic operator has rank twelve.

We can therefore write the degree-seven and degree-six equations as
\begin{equation}\label{eq:linear-elimination}
\mathcal L(H_2)=-\Mi(J_4,J_3),\qquad
\mathcal L(H_1)=-\Tri(J_4,J_3,J_2)-\det J_3.
\end{equation}
Use $(x^2,xy,xz,y^2,yz,z^2)$ for each quadratic component. The first
equation is a linear system over $\mathbb Q$ with polynomial dependence on
the parameters; Section~\ref{app:certificates} gives its compatibility
remainders. For the linear part, the selected coefficients below suffice.
Neither calculation requires division by a coefficient parameter.

\subsubsection{\texorpdfstring{The pencil $q=yz$}{The pencil q=yz}}
Source translations remove the $x$ coefficient of $u$ and the $y,z$ coefficients of $v$ in \eqref{eq:cubic-kernel}: the translation derivatives of $H_4$ correspond to $(u,v)=(2x,0),(0,z),(0,y)$. After absorbing the remaining $x$ coefficient of $v$ into the parameters, we obtain the normalized cubic
\begin{equation}\label{eq:yz-cubic}
H_3=(2x^2\ell+\alpha xq,\ \beta x^3+\gamma xq+q\ell,\ \delta x^3+\epsilon xq),
\qquad \ell=ay+bz.
\end{equation}
The first equation of \eqref{eq:linear-elimination} gives exactly
\begin{equation}\label{eq:yz-c7}
a\alpha=b\alpha=0
\end{equation}
and
\begin{equation}\label{eq:yz-quadratic}
H_2=\left(\ell^2-\gamma x\ell+Ax^2+Bq,\
\frac{6\beta-\epsilon}{4}x\ell+Cx^2+Dq,\
\frac32\delta x\ell+Ex^2+Gq\right).
\end{equation}
Write $H_1=xM+yN+zL$ for arbitrary constant vectors, and set
\[
U=\left(A+\frac{\gamma^2}{2},\
C-\frac{\gamma(6\beta+\epsilon)}8,\
E-\frac{3\delta\gamma+\epsilon^2}{4}\right)^T.
\]
Modulo \eqref{eq:yz-c7}, direct expansion gives six coefficients:
\begin{equation}\label{eq:yz-selected}
\begin{aligned}
[xy^3z^2]\Delta_6&=4(aU_1-N_1),&
[xy^2z^3]\Delta_6&=-4(bU_1-L_1),\\
[x^3y^2z]\Delta_6&=-8(aU_2-N_2),&
[x^3yz^2]\Delta_6&=8(bU_2-L_2),\\
[x^5y]\Delta_6&=4(aU_3-N_3),&
[x^5z]\Delta_6&=-4(bU_3-L_3).
\end{aligned}
\end{equation}
Hence $N=aU$, $L=bU$, and $H_1=xM+\ell U$ has Jacobian rank at most two.
This contradicts $\det\J H_1=\kappa\ne0$ and excludes the pencil $q=yz$.

\subsubsection{\texorpdfstring{The pencil $q=y^2+xz$}{The pencil q=y squared+xz}}
Here the translation derivatives correspond to $(u,v)=(2x,z),(0,2y),(0,x)$. They normalize $u=ay+bz$, $v=\rho z$. Thus, with $\ell=ay+bz$,
\begin{equation}\label{eq:parabolic-cubic}
H_3=(2x^2\ell+\alpha xq,\ q\ell+\rho x^2z+\beta x^3+\gamma xq,\ 2\rho qz+\delta x^3).
\end{equation}
The first equation of \eqref{eq:linear-elimination} gives
\begin{equation}\label{eq:parabolic-c7}
a\alpha=b\alpha=0
\end{equation}
and
\begin{equation}\label{eq:parabolic-quadratic}
\begin{split}
H_{2,1}&=\ell^2+Ax^2+Bq-\gamma x\ell+\alpha\rho xz,\\
H_{2,2}&=\rho z\ell+Cx^2+Dq+\tfrac32\beta x\ell+\gamma\rho xz,\\
H_{2,3}&=\rho^2z^2+Ex^2+Gq+\tfrac32\delta x\ell.
\end{split}
\end{equation}
Write $H_1=xM+yN+zL$ with arbitrary constant vectors.
The following two differences of coefficients give
\begin{equation}\label{eq:parabolic-two-relations}
\begin{aligned}
\bigl([x^2y^4]-[x^4z^2]\bigr)\Delta_6&=18ab\beta,\\
\bigl([x^5z]-[x^4y^2]\bigr)\Delta_6&=9ab\delta.
\end{aligned}
\end{equation}
Thus $ab\beta=ab\delta=0$. Set
\[
U=\left(A+\frac{\gamma^2}{2},\
C-\frac{3\beta\gamma}{4},\
E-\frac{3\delta\gamma}{4}\right)^T.
\]
Modulo these two relations and \eqref{eq:parabolic-c7}, three coefficients are
\begin{equation}\label{eq:parabolic-y-column}
\begin{aligned}
[x^2y^4]\Delta_6&=4(aU_1-N_1),\\
[x^4y^2]\Delta_6&=-8(aU_2-N_2),\\
[x^6]\Delta_6&=4(aU_3-N_3).
\end{aligned}
\end{equation}
Hence $N=aU$. Modulo the same relations and with this substitution,
\begin{equation}\label{eq:parabolic-selected}
\begin{aligned}
[xy^4]\Delta_5-a[x^3y^3]\Delta_6&=2a(M_1-\gamma U_1),\\
[x^3y^2]\Delta_5-2a[x^5y]\Delta_6&=-4a(M_2-\gamma U_2),\\
[x^5]\Delta_5&=2a(M_3-\gamma U_3).
\end{aligned}
\end{equation}
Therefore $aM=a\gamma U$, while $L$ remains arbitrary, and
\[
\det\J H_1=\det(M,aU,L)=\det(aM,U,L)=0.
\]
This contradicts $\det\J F(0)=\kappa$. Since no division by $a$ was used,
the argument also covers $a=0$.

\begin{proposition}\label{prop:conic-complete}
If $H_4=\nuu(p,q)$ for coprime homogeneous quadrics, a quartic Keller map is an automorphism. More precisely, a nonsquare pencil and both nonbinary square-pencil normal forms have no quartic Keller completion; the remaining binary form satisfies Proposition~\ref{prop:direction}.
\end{proposition}
\begin{proof}
Combine Proposition~\ref{prop:square-forced}, \eqref{eq:square-normalforms}, the two exclusions above, and Lemma~\ref{lem:binary-leading}.
\end{proof}

\begin{proposition}\label{prop:independent}
If the three components of $H_4$ are linearly independent, a quartic Keller map $F$ is a polynomial automorphism.
\end{proposition}
\begin{proof}
Apply Lemma~\ref{lem:factor}. For $(r,m,e)=(0,1,4)$ the leading part is binary, so Lemma~\ref{lem:binary-leading} applies. For $(1,1,3)$, take $p=x,q=y$. The linear common factor $g$ is either binary, or it can be taken as the third source coordinate. Both cases are covered by Lemma~\ref{lem:binary-leading}.

For $(2,1,2)$, the three independent binary quadrics are a basis of the binary quadratic space, so a target change gives $H_4=g(x^2,xy,y^2)$. If $g$ is independent of $z$, use Lemma~\ref{lem:binary-leading}. If $g=g_1(x,y)z+g_2(x,y)$ has $z$-degree one, its $z$ coefficient $g_1(x,y)(x^2,xy,y^2)$ has independent components, so use Lemma~\ref{lem:affine-leading}. If its $z$-degree is two, Lemma~\ref{lem:common-factor} excludes this case.

The remaining pattern $(0,2,2)$ is exactly Proposition~\ref{prop:conic-complete}. These exhaust \eqref{eq:patterns}.
\end{proof}

\subsection{A component of degree at most three: residue and boundary}\label{sec:residue}
We now treat $F=(P,Q,R)$ with
\begin{equation}\label{eq:residual-bound}
\deg P,\deg Q\le4,\qquad \deg R\le3,\qquad \det\J(P,Q,R)=\kappa\ne0.
\end{equation}
The component $R$ has no critical point. We first record how a coordinate
completion with bounded inverse degree reduces the problem to planar
invertibility.

\begin{lemma}\label{lem:completion}
Suppose, after affine source changes and affine rescaling of $R$, that $R$ is a component of a polynomial automorphism $\Theta$ with $\deg\Theta^{-1}\le5$. Then every map satisfying \eqref{eq:residual-bound} is a polynomial automorphism.
\end{lemma}
\begin{proof}
Compose with $\Theta^{-1}$ and permute the target coordinates to make the third component an affine coordinate. The new total degree is at most $4\cdot5=20$. On each fiber of the third source coordinate, the first two components form a planar Keller map of degree at most twenty, with a nonzero constant determinant. Every such map is injective by the planar invertibility theorem. Equality of the three target coordinates therefore forces equality in the source. The injectivity criterion then gives polynomial invertibility.
\end{proof}

If $\deg R\le2$, Lemma~\ref{lem:completion} already applies with inverse degree at most two. To justify the only nontrivial case, write the quadratic part as a symmetric form. If the linear part vanished on its radical, it would be in the image of the corresponding symmetric linear map and the gradient would vanish somewhere. Thus the linear part is nonzero on the radical. A linear change in that radical puts $R$ in the form $z+q(x,y)$, up to an affine constant, where $\deg q\le2$. Henceforth assume $\deg R=3$.

\subsubsection{The exact fiber form}
For each $t\in\C$ the affine surface $S_t=\{R=t\}$ is smooth. Its residue two-form $\alpha_t$ is specified by
\[
\dd R\wedge\alpha_t=\dd x\wedge\dd y\wedge\dd z,
\]
interpreted locally and then restricted to $S_t$. On $R_z\ne0$ it is $\dd x\wedge\dd y/R_z$; these expressions glue. The Keller identity gives
\begin{equation}\label{eq:residue-exact}
\dd(P|_{S_t})\wedge\dd(Q|_{S_t})=\kappa\alpha_t,
\qquad \alpha_t=\dd\bigl(\kappa^{-1}P\dd Q\bigr)|_{S_t}.
\end{equation}
Thus the residue form is the exterior derivative of the displayed regular
algebraic one-form, with the normalization fixed by $\kappa$.

\subsubsection{Boundary periods}
\begin{lemma}\label{lem:nodal-excluded}
Under \eqref{eq:residual-bound} with $\deg R=3$, the plane cubic $C=\{R_3=0\}\subset\PP^2$ is not reduced nodal. Smooth plane cubics are included in the excluded class.
\end{lemma}
\begin{proof}
Suppose $C$ is reduced nodal. Fix $t\in\C$ and let $X_t\subset\PP^3$ be the projective closure of $S_t$, with boundary $C=X_t\cap\{w=0\}$.
The surface $X_t$ is smooth along $C^{\mathrm{sm}}$: a derivative of $R_3$ tangent to the plane at infinity is nonzero there.

The form $\alpha_t$ has a simple pole along $C^{\mathrm{sm}}$, with a nowhere-vanishing holomorphic residue one-form $\eta$.
To check the pole order and normalization, use the following local calculation.
In the projective chart $(x,y,z)=(u/w,v/w,1/w)$, write $R-t=w^{-3}f$.
Then
\[
\dd x\wedge\dd y\wedge\dd z
=-w^{-4}\dd u\wedge\dd v\wedge\dd w,
\qquad
\dd R|_{f=0}=w^{-3}\dd f|_{f=0}.
\]
Consequently $\alpha_t=\dd u\wedge\dd w/(wf_v)$ where $f_v\ne0$.
Other projective charts and analytic coordinate changes preserving $w=0$ give the same expression up to a holomorphic unit. With its sign incorporated into that unit, write
\begin{equation}\label{eq:boundary-residue}
\alpha_t=e(u,v,w)\frac{\dd w\wedge\dd u}{wf_v},
\qquad
\eta=\left.\frac{e(u,v,0)\dd u}{f_v(u,v,0)}\right|_C,
\qquad e(0)\ne0.
\end{equation}
At smooth boundary points these expressions, with $u,v$ exchanged when necessary, define the asserted global residue $\eta$.

Exactness of $\alpha_t$ forces every period of $\eta$ on $C^{\mathrm{sm}}$ to vanish.
Indeed, over a closed loop $\gamma\subset C^{\mathrm{sm}}$ take the oriented normal-circle tube $T_\rho(\gamma)\subset S_t$ of small radius $\rho$.
In local coordinates $(s,w)$ on $X_t$, the form is $(\dd w/w)\wedge(a(s,w)\dd s)$ and $\eta=a(s,0)\dd s$.
Integration over the normal circle, followed by letting its radius tend to zero, gives the tube formula
\begin{equation}\label{eq:boundary-period}
0=\lim_{\rho\to0}\int_{T_\rho(\gamma)}\alpha_t
=2\pi i\int_\gamma\eta.
\end{equation}
The first equality is Stokes' theorem applied to \eqref{eq:residue-exact}. The local circle calculation patches along $\gamma$ because the residue is invariant under changes of a local defining function for $C$.

If $C$ has a node, choose analytic plane coordinates at that node and multiply the local defining equation by a unit so that $f(u,v,0)=uv$.
The numerator in \eqref{eq:boundary-residue} remains a unit.
On the punctured branch $v=0$,
\[
\eta=e(u,0,0)\frac{\dd u}{u},
\qquad
\int_{|u|=\varepsilon}\eta=2\pi i\,e(0)\ne0,
\]
contrary to \eqref{eq:boundary-period}.
If $C$ has no nodes, it is smooth and compact. Vanishing of all periods would give a global holomorphic primitive of $\eta$ on each connected component of $C$. Such a function is constant by compactness, contradicting $\eta\ne0$.
Both cases are impossible.
\end{proof}

A reduced irreducible plane cubic is smooth, nodal, or cuspidal. A reduced conic and line intersect transversely or are tangent. Three distinct lines are either nonconcurrent, giving nodes, or concurrent. A nonreduced cubic is a double line plus a line or a triple line. Consequently Lemma~\ref{lem:nodal-excluded} leaves, after linear source changes and rescaling,
\begin{equation}\label{eq:cubic-types}
R_3\in\C[x,y],\qquad R_3=x^2y+z^3,
\qquad R_3=x^2y+xz^2.
\end{equation}
The first alternative includes all repeated-line and concurrent-line configurations.

\subsubsection{A common residue obstruction}
The regular algebraic primitive in \eqref{eq:residue-exact} also gives the
following obstruction for components of a specific form.
\begin{lemma}\label{lem:period}
Let $R=x^2y+p_0(z)+xp_1(z)$, with $d=\deg p_0\ge2$ and $\deg p_1\le d$.
If the generic fibers $S_t$ are smooth and their residue forms have regular algebraic primitives, then $p_1=\lambda p_0'$ for some $\lambda\in\C$.
In particular, $R$ has a critical point.
\end{lemma}
\begin{proof}
Fix a generic $t$ and write $\alpha_t=\dd\gamma$ with $\gamma$ regular on $S_t$.
On $x\ne0$, eliminate $y$ to identify the coordinate ring with $\C[x,x^{-1},z]$. Here
\[
\alpha_t=-\frac{\dd x\wedge\dd z}{x^2}=\dd\left(\frac{\dd z}{x}\right).
\]
Every closed one-form over this Laurent ring is $\dd h+c\,\dd x/x$.
Indeed, integrate its $\dd z$ coefficient in $z$; closedness leaves a Laurent polynomial in $x$ times $\dd x$, whose only nonintegrable monomial is $x^{-1}\dd x$.
Consequently
\[
\gamma=\frac{\dd z}{x}+\dd h+c(t)\frac{\dd x}{x},
\qquad h\in\C[x,x^{-1},z].
\]

Let $z_1,\ldots,z_d$ be the distinct roots of $p_0(z)=t$, and take the local branches
\[
p_0(z_i(x))+xp_1(z_i(x))=t,\qquad z_i(0)=z_i.
\]
Pull back $\gamma$ along $(x,0,z_i(x))$. It is holomorphic at $x=0$, while the derivative of the meromorphic function $h(x,z_i(x))$ has zero residue.
Taking residues and differentiating the branch equation give
\begin{equation}\label{eq:velocities}
z_i'(0)=-c(t),\qquad p_1(z_i)=c(t)p_0'(z_i)\qquad(1\le i\le d).
\end{equation}
For small $x$, the polynomial $p_0(z)+xp_1(z)-t$ has degree $d\ge2$ in $z$. By Vi\`ete's formula, its root sum $\sum_i z_i(x)$ is independent of $t$.
Thus $-d\,c(t)=\sum_i z_i'(0)$ is constant, say $c(t)=\lambda$.
Equation~\eqref{eq:velocities}, for all generic $t$, gives $p_1=\lambda p_0'$.
At $x=0$ and any root of $p_0'$, all three partial derivatives of $R$ vanish.
\end{proof}

\subsection{Classification of the surviving cubic component}\label{sec:cubic-classification}
The boundary-period argument leaves the three types in
\eqref{eq:cubic-types}. We now examine their lower-degree terms using the
absence of critical points and the residue obstruction of
Lemma~\ref{lem:period}.
We use the following elementary elimination observation.
For $A,B,D,E\in\C[x]$ and $\lambda\in\C^*$, the equations
\[
Az+B=0,\qquad \lambda z^2+Dz+E=0
\]
have a common affine solution whenever
\begin{equation}\label{eq:quadratic-eliminant}
H(x)=\lambda B^2-DAB+EA^2
\end{equation}
is nonconstant. Indeed, choose a root $x_0$ of $H$. If $A(x_0)\ne0$, take $z=-B(x_0)/A(x_0)$. Otherwise $B(x_0)=0$, and any root of the remaining quadratic gives a solution.

\subsubsection{The cuspidal leading form}
Translations remove the quadratic coefficients of $x^2,xy,z^2$, leaving
\begin{equation}\label{eq:cusp-normal}
R=x^2y+z^3+cxz+dy^2+eyz+gx+hy+iz+j.
\end{equation}
If $d\ne0$, eliminate $y=-(x^2+ez+h)/(2d)$ from the gradient equations. The remaining equations are
\[
\begin{split}
0&=(dc-ex)z-x^3-hx+dg,\\
0&=3z^2-\frac{e^2}{2d}z+cx-\frac{e(x^2+h)}{2d}+i.
\end{split}
\]
Their eliminant \eqref{eq:quadratic-eliminant} has leading term $3x^6$, so they have an affine solution, contradicting the absence of critical points.
If $d=0,e\ne0$, solve $R_y=R_z=0$ as
\[
z=-\frac{x^2+h}{e},\qquad y=-\frac{3z^2+cx+i}{e}.
\]
Then $R_x$ becomes a polynomial in $x$ with leading term $-6x^5/e^3$, which also has a root.
Thus $d=e=0$.

If $h\ne0$, take a nonzero root of $x^2+h=0$, solve $3z^2+cx+i=0$, and then solve $2xy+cz+g=0$ for $y$. This gives a critical point. Hence $h=0$ and
\[
R=x^2y+(z^3+iz+j)+x(cz+g).
\]
Lemma~\ref{lem:period}, with $\deg p_0=3$ and $\deg p_1\le1$, now gives a critical point. This excludes the cuspidal case.

\subsubsection{A conic and its tangent line}
Translations put the component in the form
\begin{equation}\label{eq:tangent-normal}
R=x^2y+xz^2+by^2+ayz+cz^2+dx+ey+fz+j.
\end{equation}
If $b\ne0$, eliminate $y=-(x^2+az+e)/(2b)$ from the gradient equations. The remaining equations are
\[
\begin{split}
0&=\bigl(4b(x+c)-a^2\bigr)z-a(x^2+e)+2bf,\\
0&=bz^2-axz-x^3-ex+bd.
\end{split}
\]
The eliminant \eqref{eq:quadratic-eliminant} has leading term $-16b^2x^5$, so a critical point exists.
If $b=0,a\ne0$, solve $R_y=R_z=0$ as
\[
z=-\frac{x^2+e}{a},\qquad y=-\frac{2(x+c)z+f}{a}.
\]
Then $R_x$ has leading term $5x^4/a^2$ and hence has a root.
Thus $a=b=0$.
If $e\ne0$, one of the two nonzero roots of $x^2+e=0$ has $x+c\ne0$. Solve $2(x+c)z+f=0$ and then $2xy+z^2+d=0$. This is a critical point, so $e=0$. We obtain
\begin{equation}\label{eq:tangent-reduced}
R=x^2y+(x+c)z^2+dx+fz+j.
\end{equation}
If $c\ne0$, apply Lemma~\ref{lem:period} with $p_0=cz^2+fz+j$ and $p_1=z^2+d$ to obtain a critical point. Hence $c=0$.
Smoothness now forces $f\ne0$, because for $f=0$ take $x=0$ and a root of $z^2+d=0$.

The surviving scalar has an explicit coordinate completion. Set $s=xy+z^2+d$ and
\begin{equation}\label{eq:nagata-completion}
N=\left(x,\ y-\frac{2zs}{f}-\frac{xs^2}{f^2},\ z+\frac{xs}{f}\right).
\end{equation}
This map preserves $x,s$, has determinant one, and has inverse
\begin{equation}\label{eq:nagata-inverse}
N^{-1}=\left(x,\ y+\frac{2zs}{f}-\frac{xs^2}{f^2},\ z-\frac{xs}{f}\right).
\end{equation}
These statements follow by direct substitution using preservation of $s$. The inverse has degree at most five, and $R=fN_3+j$. Lemma~\ref{lem:completion} applies.

\subsubsection{A planar submersion lemma}
\begin{lemma}\label{lem:planar-submersion}
A polynomial $r\in\C[x,y]$ of degree at most three with no critical point is, up to affine source changes and affine rescaling, either of the form $Ay+h(x)$ with $A\ne0$, or $x+x^2y$.
\end{lemma}
\begin{proof}
Degrees at most two follow from the radical argument given before \eqref{eq:residue-exact}. For degree three, if the leading binary cubic is squarefree, its two gradient quadrics have no common point at infinity. B\'ezout then gives affine critical points. The remaining leading forms are $x^3$ and $x^2y$.
For leading form $x^3$, write the lower terms as $ax^2+bxy+cy^2+dx+ey+j$. If $c\ne0$, eliminate $y$ by $r_y=0$ and obtain a nonconstant quadratic equation for $x$ from $r_x=0$. If $c=0,b\ne0$, solve $r_y=0$ for $x$ and then $r_x=0$ for $y$. Both give critical points. Therefore $b=c=0$, and smoothness forces $e\ne0$, giving the graph form.
For leading form $x^2y$, translations leave
\[
r=x^2y+Cy^2+Dx+Ey+J.
\]
If $C\ne0$, solve $r_y=0$ for $y$; substitution in $r_x=0$ gives the nonconstant cubic equation $-x(x^2+E)/C+D=0$. If $C=0,E\ne0$, choose a nonzero root of $x^2+E=0$ and solve $2xy+D=0$. Hence $C=E=0$, and absence of a critical point is exactly $D\ne0$. Affine rescaling gives $x+x^2y$.
\end{proof}

\subsubsection{Binary cubic leading parts}
When $R_3\in\C[x,y]$, the component has the form
\begin{equation}\label{eq:binary-scalar}
R=az^2+\ell(x,y)z+p(x,y),\qquad \deg\ell\le1,\quad\deg p\le3.
\end{equation}
If $a\ne0$, an affine completion of the square gives $az^2+\widetilde p(x,y)$. The planar polynomial $\widetilde p$ has no critical point. Lemma~\ref{lem:planar-submersion} applies. Its graph case gives a coordinate completion with inverse degree at most three. The other case, after rescaling, is
\[
R=x+x^2y+az^2,\qquad a\ne0.
\]
Lemma~\ref{lem:period}, with $p_0=az^2$ and $p_1=1$, contradicts the absence of critical points.

If $a=0$ and $\ell$ is a nonzero constant, \eqref{eq:binary-scalar} is a graph coordinate, again with inverse degree at most three. If $a=0$ and $\ell$ is nonconstant, take $\ell=x$. For $R=xz+p(x,y)$, absence of critical points requires $p_y(0,y)$ to be a nonzero constant: otherwise choose a root at $x=0$ and then solve $R_x=0$ for $z$. Hence
\[
p=Ay+B+xh(x,y),\qquad A\ne0,\quad\deg h\le2.
\]
Put $W=z+h(x,y)$. The completion $(x,R,W)$ has inverse, in coordinates $(X,T,W)$,
\begin{equation}\label{eq:binary-completion}
x=X,\qquad y=\frac{T-XW-B}{A},\qquad
z=W-h\left(X,\frac{T-XW-B}{A}\right),
\end{equation}
of degree at most four. Finally, $a=0$, $\ell=0$ is precisely the planar case of Lemma~\ref{lem:planar-submersion}.

These arguments establish the following reduction under the Keller
hypotheses.
\begin{proposition}\label{prop:cubic-reduction}
Under \eqref{eq:residual-bound}, either $R$ has a polynomial coordinate completion whose inverse has degree at most five, after affine changes and rescaling, or $R$ is affine-equivalent to $x+x^2y$.
\end{proposition}
\begin{proof}
For degrees at most two use the radical argument. For degree three, Lemma~\ref{lem:nodal-excluded} gives \eqref{eq:cubic-types}. The cuspidal case has been excluded; the tangent case gives \eqref{eq:nagata-completion}, and the binary case gives \eqref{eq:binary-completion}, graph coordinates, or the stated exception.
\end{proof}

\subsection{\texorpdfstring{The exceptional scalar $x+x^2y$}{The exceptional scalar x+x squared y}}\label{sec:exceptional-scalar}
It remains to exclude the exceptional scalar in
Proposition~\ref{prop:cubic-reduction}. We express the Keller determinant
through a derivation tangent to its fibers. A formal residue obstruction
then restricts the leading coefficients of a possible completion.
\begin{proposition}\label{prop:exceptional-scalar}
There are no $P,Q\in\C[x,y,z]$ with $\deg P,\deg Q\le4$ and
\[
\det\J(P,Q,x+x^2y)\in\C^*.
\]
\end{proposition}
\begin{proof}
Set $B=\C[x,y]$, $r=x+x^2y$, and
\[
D=-x^2\partial_x+(1+2xy)\partial_y.
\]
For $P,Q\in B[z]$ define $[P,Q]_D=D(P)Q_z-P_zD(Q)$. Direct expansion gives
\begin{equation}\label{eq:exceptional-bracket}
\det\J(P,Q,r)=[P,Q]_D.
\end{equation}
Assume $[P,Q]_D=\kappa\ne0$.

\subsubsection*{No polynomial slice}
With $u=1/x$ one has
\[
x=u^{-1},\qquad y=ru^2-u,\qquad
\C(x,y)=\C(r,u),\qquad D=\partial_u.
\]
Thus
\begin{equation}\label{eq:constant-field}
\Ker_{\C(x,y)}D=\C(r),\qquad D(1/x)=1.
\end{equation}
If $Dg=c\in\C^*$ with $g\in B$, then $g=c/x+h(r)$ for $h\in\C(t)$.
Write $h=A/B_0$ in reduced form.
Since $h(r)=(xg-c)/x$, the divisibility observation after \eqref{eq:intersection} gives $B_0(r)\mid x$.
But $\deg B_0(r)=3\deg B_0$, so $B_0$ is constant.
This would make $c/x=g-h(r)$ polynomial, a contradiction. Therefore
\begin{equation}\label{eq:no-slice}
Dg\notin\C^*\quad\hbox{for every }g\in B.
\end{equation}

\subsubsection*{A coefficient at infinity}
Suppose $P$ has positive $z$-degree $m$ and nonzero constant leading coefficient. Rescale it to be monic and form the binomial series
\[
S=P^{-1/m}=z^{-1}(1+O(z^{-1}))\in B((z^{-1})).
\]
Write $\operatorname{res}_z$ for the coefficient of $z^{-1}$, so $\operatorname{res}_zS=1$ and the residue of every formal $z$-derivative is zero.
Extend $D$ coefficientwise and set $g=\operatorname{res}_z(QP_zS)\in B$.
The chain rule gives
\begin{equation}\label{eq:formal-slice}
D(QP_zS)-\partial_z(QD(P)S)=-\kappa S.
\end{equation}
Taking residues gives $Dg=-\kappa$, contrary to \eqref{eq:no-slice}.
The argument applies to either component. A component independent of $z$ is also impossible: the bracket equation would make its $D$-derivative a nonzero constant.
Thus both $z$-degrees are positive and their leading coefficients are nonconstant.

\subsubsection*{No affine pair, in any total degree}
If $P=pz+a$ and $Q=qz+b$, Lemma~\ref{lem:affine-reduction} gives
\[
D(qa-pb)=\{pb-qa,r\}=\kappa,
\]
contrary to \eqref{eq:no-slice}.
This argument requires no total degree bound and allows either $p$ or $q$
to vanish.

\subsubsection*{Leading coefficients and squarefree fibers}
Write
\[
m=\deg_zP,\quad n=\deg_zQ,\quad
p=\operatorname{lc}_zP,\quad q=\operatorname{lc}_zQ.
\]
The highest $z$ coefficient of the bracket gives
\begin{equation}\label{eq:lead-ratio}
nqD(p)-mpD(q)=0,
\qquad \frac{p^n}{q^m}\in\C(r).
\end{equation}
Interchange $P,Q$ if needed so that $m\ge n$.
The preceding obstructions give $2\le m\le3$ and
\[
1\le\deg p\le4-m\le2,\qquad \deg q\le4-n\le3.
\]

We claim that the ratio in \eqref{eq:lead-ratio} is constant.
Every polynomial $r-c$, $c\in\C$, is squarefree of degree three: its partial derivatives $1+2xy$ and $x^2$ have no common zero.
Write
\[
\frac{p^n}{q^m}=\frac{A(r)}{B_0(r)},\qquad
A,B_0\in\C[t],\quad\gcd(A,B_0)=1.
\]
B\'ezout gives $\gcd(A(r),B_0(r))=1$, so $A(r)\mid p^n$ and $B_0(r)\mid q^m$.
If $A$ has a root $c$, squarefreeness gives $r-c\mid p$, contrary to $\deg p\le2$.
If $B_0$ has a root $c$, it similarly gives $r-c\mid q$; the degree bound forces $q=b(r-c)$ for $b\in\C^*$.
Then $D(q)=0$, and \eqref{eq:lead-ratio} gives $D(p)=0$.
By \eqref{eq:constant-field} and \eqref{eq:intersection}, $p\in\C[r]$. Its degree bound would make it constant, again a contradiction.
Thus $A,B_0$ are constant, proving the claim.

If $m=n$, the claim gives $p/q\in\C^*$.
A constant target shear lowers $m+n$ without increasing total degree.
Repeat this reduction until $m>n$, applying the preceding obstructions at each step to exclude zero $z$-degree, constant leading coefficients, and affine pairs.
Since $m\le3$, the positive integers $m,n$ are then coprime.
Unique factorization in $p^n=cq^m$ gives
\[
p=a\ell^m,\qquad q=b\ell^n,\qquad
a,b\in\C^*,\quad \ell\in B\setminus\C.
\]
The inequality
\[
m\le\deg p\le4-m
\]
forces $m=2$, and hence $n=1$ and $p=\lambda q^2$ for $\lambda\in\C^*$.
Now $U=P-\lambda Q^2$ and $Q$ are both affine in $z$, with $[U,Q]_D=\kappa$.
The affine-pair obstruction, which requires no total degree bound, gives
the final contradiction even if the total degree of $U$ increases.
\end{proof}

\subsection{Proof of invertibility}\label{sec:assembly}
The independent-leading-part case was settled in
Proposition~\ref{prop:independent}. We combine the remaining scalar
reductions and then complete the proof of Theorem~\ref{thm:main}.
\begin{proposition}\label{prop:low-component}
A map satisfying \eqref{eq:residual-bound} is a polynomial automorphism.
\end{proposition}
\begin{proof}
Proposition~\ref{prop:cubic-reduction} gives either a coordinate completion of inverse degree at most five or the exceptional scalar, after affine changes. In the first case apply Lemma~\ref{lem:completion}; in the second apply Proposition~\ref{prop:exceptional-scalar}. Affine changes preserve the total degree bounds.
\end{proof}

\begin{proof}[Proof of Theorem~\ref{thm:main}]
Subtract the constant term of $F$ and let $H_4$ be its quartic homogeneous part. If the three components of $H_4$ are linearly independent, Proposition~\ref{prop:independent} applies. Otherwise there is a nonzero constant linear form on the target annihilating $H_4$. Complete it to an invertible constant target transformation and write the resulting map as $(P,Q,R)$ with $\deg R\le3$. The component $R$ cannot be constant because the Jacobian is nonzero. The other two degrees are still at most four. Proposition~\ref{prop:low-component} applies. This includes $H_4=0$ and therefore degrees smaller than four. Reversing the affine transformations proves the assertion.
\end{proof}

\begin{corollary}\label{cor:charzero}
The conclusion of Theorem~\ref{thm:main} holds over every characteristic-zero field.
\end{corollary}
\begin{proof}
The coefficients of the map lie in a finitely generated extension $K_0$ of $\mathbb Q$, which embeds in $\C$. The complex result gives a polynomial inverse after that embedding. Translate and normalize at the origin over $K_0$. The formal inverse has coefficients in $K_0$, by recursive inversion of its invertible linear part. Uniqueness identifies it with the complex polynomial inverse, so it has finite support and belongs to $K_0[\XX]^3$. The two inverse identities then hold over $K_0$ and after any extension of scalars.
\end{proof}

\subsection{Exact coefficient certificates}\label{app:certificates}
We record the coefficient operators and compatibility remainders used in
Section~\ref{sec:square}. The quadratic equations are linear systems over
$\mathbb Q$ with right-hand sides polynomial in the parameters. The
linear-part obstructions use the selected coefficients displayed there.
All calculations avoid division by coefficient parameters.

\subsubsection{The operators defining the matrices}
For $q=yz$ or $q=y^2+xz$, define $n=(q^2,-2x^2q,x^4)^T$ and $\mathcal L(U)=4x n^T\calE U$ as in \eqref{eq:coefficient-operator}. This formula explicitly defines every entry of the matrices used in \eqref{eq:linear-elimination}; rows are the coefficients of degree-seven or degree-six monomials, ordered lexicographically with $x>y>z$. Columns are the coefficients in each component, with quadratic order
\[
(x^2,xy,xz,y^2,yz,z^2)
\]
and linear order $(x,y,z)$. The elementary actions are
\begin{center}
\begin{tabular}{c|rrrrrr}
& $x^2$&$xy$&$xz$&$y^2$&$yz$&$z^2$\\\hline
$z\partial_z-y\partial_y$&$0$&$-xy$&$xz$&$-2y^2$&$0$&$2z^2$\\
$-x\partial_y+2y\partial_z$&$0$&$-x^2$&$2xy$&$-2xy$&$2y^2-xz$&$4yz$
\end{tabular}
\end{center}
The kernel \eqref{eq:quadratic-kernel} shows that the quadratic operator
has rank twelve for all parameter values.

\subsubsection{Quadratic elimination}
Write $H_2=(\sum_{j=0}^5c_jm_j,\sum_{j=0}^5c_{6+j}m_j,\sum_{j=0}^5c_{12+j}m_j)$ in the listed quadratic monomial basis. For $q=yz$, rational Gaussian elimination of the first equation in \eqref{eq:linear-elimination} leaves the free coefficients
\[
c_0=A,\ c_4=B+2ab,\ c_6=C,\ c_{10}=D,\ c_{12}=E,\ c_{16}=G.
\]
The two compatibility remainders are
\[
-2a\alpha,\qquad 2b\alpha.
\]
Solving the remaining twelve equations gives exactly \eqref{eq:yz-quadratic}.
For $q=y^2+xz$, take
\[
c_0=A,\ c_3=B+a^2,\ c_6=C,\ c_9=D,\ c_{12}=E,\ c_{15}=G.
\]
The compatibility remainders, with repetitions retained, are
\[
-6a\alpha,\ -12a\alpha,\ 4b\alpha,\ -6a\alpha,\ 8b\alpha,\ 4b\alpha.
\]
Modulo their ideal, the remaining twelve equations give \eqref{eq:parabolic-quadratic}. These lists also show directly why the necessary parameter ideal is $(a\alpha,b\alpha)$ in both cases.

\section{Tameness of quartic automorphisms}\label{sec:tameness}
We now prove tameness under an independent automorphism assumption.
The argument combines a directional criterion with the classification of
low-degree coordinates and local tameness of plane polynomial families.
The leading-part calculations in Section~\ref{sec:invertibility} will
also be used, without repeating their coefficient computations.

\begin{theorem}[Tameness]\label{thm:tame-aut}
Every complex polynomial automorphism of $\A^3$ of total degree at most
four is tame.
\end{theorem}

\subsection{The directional criterion for tameness}\label{t:direction}
We first strengthen the directional invertibility criterion of
Proposition~\ref{prop:direction} under the assumption that the map is an
automorphism. The resulting tameness criterion requires no degree bound.

\begin{proposition}\label{t:directional}
Let $F\in\Aut_{\C}(\A^3)$, with no degree assumption.
If $\partial_e^2F$ is constant for some nonzero constant vector $e$, then $F$ is tame.
\end{proposition}
\begin{proof}
After a linear source change, take $\partial_e=\partial_z$ and write
\[
 F=A(x,y)+zB(x,y)+z^2E,
 \qquad E\in\C^3.
\]
Every nonzero linear combination $\ell(F)$ of the components is a
coordinate, so $\ell(F)-d$ is irreducible for each $d\in\C$.

First suppose $E=0$. For $\tau_t(x,y,z)=(x,y,z+t)$, conjugation gives
\[
 F\circ\tau_t\circ F^{-1}(Y)=Y+tH(Y),\qquad H=B\circ F^{-1}.
\]
Thus $Y+H(Y)$ has inverse $Y-H(Y)$: it is a quasi-translation.
Moreover, $H$ has no zero because $B=F_z$ is a column of an invertible
Jacobian matrix. The classification in \cite[Corollary 5.8]{QT} gives,
after a linear target change, $H=g\cdot(p(W),q(W),0)$ in target coordinates
$(U,V,W)$, where $p,q\in\C[W]$ and $g\in\C[U,V,W]$.
Since $H$ has no zero, $g$ is a nonzero constant. Absorb it into $p,q$;
then $p,q$ have no common root. Pulling back by $F$ gives
\[
 F=(a+zp,b+zq,r),\qquad a,b,r\in\C[x,y],\quad p,q\in\C[r].
\]
Choose $\sigma,\tau\in\C[r]$ such that $\sigma p+\tau q=1$. The matrix
\[
 \begin{pmatrix}\sigma&\tau\\-q&p\end{pmatrix}\in\mathrm{SL}_2(\C[r])
\]
is elementary over the Euclidean ring $\C[r]$. Its tame target action sends $F$ to
$(z+\alpha(x,y),S(x,y),r(x,y))$. Polynomial invertibility implies that $(S,r)$
is a plane automorphism: evaluate the first two components of the inverse with
the first target coordinate set to zero. The Jung--van der Kulk theorem completes this case.

Now suppose $E\ne0$. Normalize
\[
 F=(a+zp,b+zq,c+zr+z^2).
\]
Again $p,q$ have no common zero, since at such a zero $z=-r/2$ would give $F_z=0$.
Complete the square by an elementary source change, retaining $p,q$ and
relabeling $a,b,c$. The map is now $(a+zp,b+zq,c+z^2)$.
With $\br fg=f_xg_y-f_yg_x$, the two highest determinant coefficients give
\begin{equation}\label{t:Eeqs}
 \br pq=0,\qquad \br pb+\br aq=0.
\end{equation}
Suppose no constant combination of $p,q$ is a nonzero constant.
Lemma~\ref{lem:differential-constants} and
Proposition~\ref{prop:closed-generators}(iii) give a closed common
generator $s$, so $p=P(s)$, $q=Q(s)$.
The polynomial kernel identity in Lemma~\ref{lem:differential-constants} gives
$\ker\br s{\,\cdot\,}=\C[s]$.
The polynomials $P,Q$ are coprime; a constant target change gives
$0<m=\deg P<n=\deg Q$.
Put $g=\gcd(P',Q')$ and $(P_1,Q_1)=(P',Q')/g$.
The second equation in \eqref{t:Eeqs} now gives
\[
 \br s{P_1(s)b-Q_1(s)a}=0,
 \qquad P_1(s)b-Q_1(s)a=h(s).
\]
Since
\[
 \deg(PQ_1-QP_1)=m+n-1-\deg g\ge n>0,
\]
choose a root $c_0$ of $PQ_1-QP_1$.
The vectors $(P(c_0),Q(c_0))$ and $(P_1(c_0),Q_1(c_0))$ are nonzero
and proportional. Modulo $s-c_0$, the displayed identity therefore makes
$Q(c_0)F_1-P(c_0)F_2$ constant.
The irreducibility observation forces this coordinate minus its constant
to equal $\lambda(s-c_0)$. Its $z$ coefficient is consequently zero:
$Q(c_0)P(s)-P(c_0)Q(s)=0$, contradicting $m<n$.

Thus a constant target change gives a component $t=z+a_0(x,y)$.
Straightening this component puts the transformed automorphism in the form
\[
 \bigl(A(u,v)+tB(u,v)+t^2E',\ t\bigr),\qquad E'\in\C^2.
\]
Subtract $t^2E'$ from the first two target components by elementary shears.
The resulting automorphism is affine in $t$, so the already proved $E=0$
case applies; it requires no degree bound.
\end{proof}

\subsection{Small coordinate components and polynomial families}\label{t:small}
Suppose the quartic leading components span a space of dimension at most
two. A linear target change then produces a component of degree at most
three. We combine the low-degree coordinate classification with local
tameness over $\C[t]$, beginning with an exclusion needed in both parts
of the argument.

\subsubsection{An elementary exclusion of the Nagata component}

\begin{lemma}\label{t:nagata}
The coordinate $r_d=z+x(xy+z^2+d)$ cannot be completed to an automorphism
$(P,Q,r_d)$ with $\deg P,\deg Q\le4$.
\end{lemma}
\begin{proof}
Set $s=xy+z^2+d$, $q=y-2zs-xs^2$, and $r=z+xs$.
The automorphism $(x,q,r)$ preserves $s$, and $s=xq+r^2+d$.
Work on the generic fiber $r=c$ over $K=\C(c)$, and put $D=c^2+d$.
Its inverse coordinates are
\begin{equation}\label{t:fiber}
 z=c-x(xq+D),\qquad y=q+2c(xq+D)-x(xq+D)^2.
\end{equation}
The inverse in \eqref{t:fiber} is affine on each coordinate axis, so both
axis restrictions of an ambient quartic have degree at most four.
The restrictions $\bar P,\bar Q$ are plane coordinates over $K$.
We use the following consequence of the Jung--van der Kulk theorem:
over an algebraically closed field, the highest homogeneous part of a
plane coordinate is a nonzero scalar multiple of a power of a linear
form. Indeed, in a reduced affine--triangular decomposition, each
nonlinear degree-increasing step replaces the highest form of one
component by a scalar multiple of a power of the other. Starting with
an affine pair gives the assertion; a final affine combination either
retains that highest form or recovers the preceding component.
Apply this observation after extending scalars from $K$ to $\overline K$,
which preserves both the coordinate property and degree. If the highest
form of a restricted coordinate of degree $m$ is
$\lambda(\alpha x+\beta q)^m$, then $\alpha,\beta$ are not both zero.
Its degree is therefore attained on the $x$-axis or the $q$-axis.
The bounds on the two axis restrictions consequently give
\begin{equation}\label{t:fiberdegree}
 \deg\bar P,\deg\bar Q\le4.
\end{equation}

Use the Laurent curve
\[
 (x,y,z)=(-t,D/t,t^{-1}-c),\qquad r=c,\qquad
 q=-t^{-3}+2ct^{-2}+Dt^{-1}.
\]
Every ambient quartic has pole order at most four at $t=0$.
The leading Laurent term of $x^iq^j$ is $(-1)^{i+j}t^{i-3j}$.
Let $u$ be either restricted coordinate and $m=\deg u\le4$.
If $m\ge2$, the monomial $q^m$ has the unique greatest pole order $3m>4$,
so its coefficient vanishes. The pure-power leading form of $u$ must
therefore be a multiple of $x^m$. The remaining terms have degree at most
three, where the pole orders $3j-i$ are pairwise distinct; any term with
$j\ge2$ has pole order at least five and hence vanishes.
Thus, also when $m=1$, we have $u=A(x)q+h(x)$ over $K$.
If $A$ is nonconstant, choose a root $a$ over $\overline K$.
Then $x-a$ divides the irreducible coordinate $u-h(a)$, forcing
$u-h(a)\in\overline K^*(x-a)$, contrary to $A\ne0$.
Hence $A$ is constant.

The pullbacks of $P,Q$ lie in $\C[x,q,r]$, so each has the form
\[
 A(r)q+h(x,r),\qquad A\in\C[r],\quad h\in\C[x,r].
\]
Set $S=xy+z^2$. The leading forms of $r,q$ are $xS,-xS^2$.
Since $x,S$ are algebraically independent, the leading form of $h(x,r)$
lies in $\C[x,xS]$, whose monomials $x^iS^j$ have $i\ge j$.
For $A\ne0$ of degree $n$, the leading form of $A(r)q$ is proportional to
$x^{n+1}S^{n+2}$ and cannot cancel. Thus
$\deg(A(r)q+h(x,r))\ge3n+5>4$, a contradiction.
Both $P,Q$ therefore lie in $\C[x,r]$, contradicting invertibility of $(P,Q,r)$.
\end{proof}

\subsubsection{The cases \texorpdfstring{$\rho\le1$}{rho at most one}}\label{t:smallpair}
The cases with at most one independent quartic leading component follow
from the degree-three classification and the preceding exclusion.
If $\rho=0$, tameness follows from \cite[Corollary~3.9.5]{BVS}.
If $\rho=1$, choose two components of degree at most three. The inverse image
of each affine line under this pair is $\A^2$. The characteristic-zero cases
(4)--(7) of \cite[Theorem 3]{BVS} apply. The first three have tame completions;
the last contains the coordinate of Lemma~\ref{t:nagata} and is excluded.
For the first three, the completions are triangular or are obtained from
\[
 (x+yz+za(x,z),\ y+a(x,z)+r(z),\ z)
\]
by projection. The displayed automorphism is a composition of three shears.
Any other completion of the same pair differs by an affine automorphism of
the remaining variable over the polynomial ring in that pair, and is also tame.

\subsubsection{Plane automorphisms over \texorpdfstring{$\C[t]$}{C[t]}}
For the remaining coordinate forms, we will use the following local
criteria and their passage to global tameness.

\begin{lemma}\label{t:families}
Let $\Phi\in\GA_2(\C[t])$.
\begin{enumerate}[label=(\roman*)]
\item At a fiber where the maximum plane degree equals the generic degree,
$\Phi$ is tame over the corresponding local ring.
\item Over a DVR, every plane automorphism of generic degree at most three is tame.
\item If $\Phi$ is tame at every closed point of $\A^1$, it is tame over $\C[t]$.
\end{enumerate}
\end{lemma}
\begin{proof}
For (i), let $R=\C[t]_{(t-c)}$, let $K=\C(t)$, and write a bar
for specialization at $t=c$. A polynomial over $R$ is called
\emph{primitive} if at least one coefficient is a unit. Thus the
highest homogeneous part of a component is primitive precisely when
its degree is preserved by specialization. Any finite sequence of
transformations over $R$ is defined on a sufficiently small connected
open neighborhood of $c$.

Here \emph{multidegree} denotes the ordered sequence of degrees of the
nonlinear triangular factors in a reduced Jung--van der Kulk
decomposition, rather than the pair of component degrees
$(\deg p,\deg q)$.
We first recall two consequences of this normal form over a field.
For an automorphism with unequal component
degrees $1<d<D$, one has $D=jd$ with $j\ge2$, and the higher-degree
component can be reduced to degree less than $d$ by subtracting a
polynomial in the lower-degree component. Moreover, the degree of this
last triangular factor is $j=D/d$; it is the entry of the multidegree
removed first by this degree reduction. Consequently, for two
unequal-degree pairs with the same multidegree and the same maximum
degree $D$, the smaller degree is the same. These assertions follow
by reading the final triangular factor of a reduced decomposition.

Write $\Phi=(p,q)$ and induct on $\deg p+\deg q$, maintaining the
hypothesis that the generic and special maximum degrees agree.
If either component is linear, the straightening argument in (ii) below
already proves tameness. We may therefore assume that both generic
component degrees exceed one. After shrinking around $c$, the degree
is constant on a connected open neighborhood. By
\cite[Corollary~3]{F}, the multidegree is therefore constant in this
family, so the generic and special fibers have the same multidegree.
If both generic degrees equal $D>1$,
choose the component whose highest form is primitive and call it $q$;
such a component exists because the special maximum degree is $D$.
Over $K$, the highest forms satisfy $p_D=\lambda q_D$.
A unit coefficient of $q_D$ shows that $\lambda\in R$.
The elementary change $(p,q)\mapsto(p-\lambda q,q)$ reduces the sum
of the generic degrees, while $q$ retains degree $D$ in both fibers.
The induction hypothesis applies to this new pair.

It remains to consider unequal generic degrees $d=\deg q<D=\deg p$.
The specialization of $p$ has degree $D$, and that of $q$ has degree
at most $d<D$. It cannot become linear: the special multidegree would
then have the single entry $D$, whereas $d>1$ gives at least two
entries in the generic multidegree. Constancy of multidegree and the
preceding normal-form observation therefore force $\deg\bar q=d$.
Hence $q_d$ is primitive.
Over $K$, subtract the polynomial in $q$ that reduces $p$ below
degree $d$. Each coefficient in this polynomial belongs to $R$:
at a step with current degree $jd\ge d$, the highest-form identity is
$p_{jd}=\lambda q_d^j$. Since $q_d$ is primitive, its reduction is
nonzero, and so is the reduction of $q_d^j$. Thus $q_d^j$ has a unit
coefficient, and $\lambda\in R$. Subtracting $\lambda q^j$ preserves
integrality of all coefficients. Repeating, including a linear term
when needed, gives $p'=p-A(q)$ with $A\in R[T]$ and $\deg p'<d$.
The pair $(p',q)$ has generic and special maximum degree $d$, since
$q_d$ remains primitive. Its degree sum is smaller, so induction
applies. The process terminates at an affine automorphism over $R$.
This proves (i).

For (ii), let $R$ be a DVR with fraction field $K$. If a component is
linear, one of its linear coefficients is a unit: otherwise its
specialization would be constant, contrary to invertibility of the
specialized map. An affine change over $R$ therefore makes this
component a source variable. Over $K$, the other component is then
linear in the remaining variable; its coefficients lie in $R$ and its
linear coefficient is a unit, since the Jacobian determinant is a unit.
Thus the resulting automorphism is triangular over $R$.
If neither component is linear, their degrees are $(2,2)$ or $(3,3)$;
the pair $(2,3)$ is excluded by the plane degree-divisibility theorem.
Their highest forms are proportional over $K$. For the proportionality
coefficient, either it or its reciprocal lies in $R$, so an elementary
linear change over $R$ lowers one component degree. The resulting
unequal degrees force a linear component, again by divisibility.
This reduces to the preceding case.

Finally, (iii) follows from \cite[Theorem 3]{BDFM}, since every invertible
module over $\C[t]$ is free. The generic localization is a field and
hence is tame by the Jung--van der Kulk theorem, so checking the closed
points suffices.
\end{proof}

\subsubsection{Components of the form \texorpdfstring{$x+C(y,z)$}{x+C(y,z)}}
Assume $\rho=2$ and take $F=(P,Q,R)$ with $\deg R\le3$.
A quadratic coordinate has no critical point. Its linear part is nonzero on
the radical of its quadratic part, so affine changes put it in one of the forms
\[
 R=x,\qquad R=x+y^2,\qquad\text{or}\qquad R=x+yz.
\]
For a cubic coordinate, \cite[Theorem 3, I]{BVS} and Lemma~\ref{t:nagata}
leave, up to affine equivalence, precisely the forms
\[
 R=x+C_2(y,z)+C_3(y,z),\qquad R=xy+yC_2(y,z)+z,
\]
where $C_i$ is homogeneous of degree $i$ and, in the second form,
$C_2\notin\C[y]$. We treat the second form in Section~\ref{t:noncube}.
All the other forms have $R=x+C(y,z)$. Straighten them by $t=R$;
the remaining pair lies in $\GA_2(\C[t])$.

We isolate the degree estimate that controls specialization of these
families. Here relative degree means degree in $y,z$, with $t$ treated
as a coefficient.
\begin{lemma}[Degree drop under a graph substitution]\label{lem:graph-degree-drop}
Let $R=x+C(y,z)$ have degree $m\in\{1,2,3\}$, and let
$A\in\C[x,y,z]$ satisfy $\deg A\le4$. Suppose
$a(t,y,z)=A(t-C,y,z)$ has generic relative degree $D$ and
$\deg_{y,z}a(c,y,z)<D$ for some $c\in\C$. Then there are polynomials
$B,A_0\in\C[x,y,z]$ such that
\[
 A=(R-c)B+A_0,\qquad \deg B\le4-m,
 \qquad \deg_{y,z}A_0(t-C,y,z)<D,
\]
and
\begin{equation}\label{t:dropbound}
 D=\deg_{y,z}B(t-C,y,z)\le m(4-m).
\end{equation}
\end{lemma}
\begin{proof}
Give $x,y,z$ weights $m,1,1$, and let $C_m$ be the degree-$m$ part
of $C$ (possibly zero when $m=1$). Start with remainder $A$ and
quotient zero. If the highest weighted form $W$ of the current
remainder is divisible by $x+C_m$, write $W=(x+C_m)U$ and subtract
$(R-c)U$ from that remainder, adding $U$ to the quotient.
The polynomial $U$ is weighted homogeneous. Also, $W$ has ordinary
degree at most four, and $x+C_m$ has ordinary degree $m$. Additivity
of ordinary degree under multiplication gives
$\deg U=\deg W-m\le4-m$. Thus the subtraction preserves the
ordinary degree bound four. It strictly decreases the highest weight,
because all terms of $R-c-(x+C_m)$ have weight less than $m$.
The process therefore terminates, giving
$A=(R-c)B+A_0$ with $\deg B\le4-m$ and $\deg A_0\le4$.

If $A_0\ne0$, its highest weighted form $W_0$ is not divisible by
$x+C_m$. The kernel of substitution $x=-C_m$ is the principal ideal
$(x+C_m)$, so $W_0(-C_m,y,z)\ne0$. If $w$ is its weight, this
substitution is homogeneous of ordinary degree $w$ in $y,z$.
It is the highest homogeneous part of both $A_0(t-C,y,z)$ and
$A_0(c-C,y,z)$: the terms involving $t$, the constant $c$, or a
lower-degree part of $C$ have smaller relative degree. Hence these
two polynomials have the same relative degree. Since
$a(c,y,z)=A_0(c-C,y,z)$, that degree is less than $D$ by hypothesis.
The conclusion is immediate as well if $A_0=0$.

Substituting $x=t-C$ in the division identity now gives
$a(t,y,z)=(t-c)B(t-C,y,z)+A_0(t-C,y,z)$. The remainder has relative
degree less than $D$, so $B(t-C,y,z)$ has relative degree exactly
$D$. Finally, substitution increases ordinary degree by at most a
factor of $m$, giving \eqref{t:dropbound}.
\end{proof}

Apply Lemma~\ref{lem:graph-degree-drop} to each component whose relative
degree drops at the fiber under consideration.
For $m=1$ or $3$, every drop has generic degree at most three, so
Lemma~\ref{t:families}(i)--(iii) proves tameness. This treats every cubic
$C$, independently of the factorization of $C_3$.
It remains to consider $m=2$, with $C=y^2$ or $yz$.
For homogeneous quadratic $C$, a degree-four drop comes only from the $x^2$ term of $B$,
and its leading form is
\begin{equation}\label{t:quadtop}
 b(t-c)C^2,\qquad b\in\C^*.
\end{equation}

At a fiber where the maximum degree drops, take $c=0$.
Lemma~\ref{t:families} handles generic degrees
at most three. A degree-four dropping coordinate has leading form $btC^2$,
which must be a scalar multiple of a fourth power. Thus $C=yz$ is impossible,
and we take $C=y^2$. If both generic degrees are four, a constant target shear
cancels one leading part. Plane degree divisibility then gives degrees
$(4,1)$ or $(4,2)$; write the smaller component as $q$.
In the latter case $p_4=\lambda q_2^2$, so $q_2\in\C(t)C$; in the former
$q_2=0$.

Set $U=q_t(0,y,z)$ and $V=q(0,y,z)$, both of degree at most two.
Taylor's formula gives $R^2\mid Q-RU-V$; the degree bound therefore yields
$Q=aR^2+RU+V$ for some $a\in\C$.
Writing $U_2,V_2$ for their homogeneous quadratic parts, we have
$tU_2+V_2=q_2\in\C(t)C$, so $U_2\in\Span_{\C}\{C\}$.
Consequently the ambient quartic part satisfies
\begin{equation}\label{t:quadraticreduction}
 Q_4=aC^2+CU_2\in\Span_{\C}\{C^2\}=\Span_{\C}\{R_2^2\}.
\end{equation}
Subtracting a constant multiple of $R^2$ from $Q$ cancels $Q_4$.
The transformed automorphism still has total degree at most four and now
has two components of degree at most three,
so the argument in Section~\ref{t:smallpair} proves tameness. If no degree-four drop occurs, every fiber is
covered by Lemma~\ref{t:families}(i)--(ii), and part (iii) gives global tameness.

\subsection{The remaining cubic coordinate form}\label{t:noncube}
The preceding argument treats components of the form $x+C(y,z)$.
We now handle the other cubic normal form, using invariants of its leading
binary cubic to obtain a constant second directional derivative.

Write $R=xy+yC_2(y,z)+z$, where $C_2\notin\C[y]$.
Thus $R_{xx}=0$ and $h=R_3=yC_2$ is a noncube binary cubic.
Set $D=h_z\partial_y-h_y\partial_z$.
For a nonzero homogeneous rational invariant $u$ of degree $d$, Euler's
identity gives $3h\,\mathrm du=d\,u\,\mathrm dh$, hence $u^3/h^d\in\C^*$.
A noncube binary cubic has a simple linear factor. Its valuation forces
$3\mid d$, and then $u\in\C^*h^{d/3}$. Taking homogeneous components gives
\begin{equation}\label{t:hinvariants}
 \ker_{\C[y,z]}D=\C[h].
\end{equation}
The independent quartics $f=P_4,g=Q_4$ satisfy
\begin{equation}\label{t:pairleading}
 f_xDg-g_xDf=0.
\end{equation}

Equal positive $x$-degrees in \eqref{t:pairleading} can be reduced by constant
target shears: the leading-coefficient ratio is a degree-zero homogeneous
$D$-invariant, hence constant. If exactly one $x$-degree is zero, its
nonzero binary quartic is a $D$-invariant, impossible.
If the unequal degrees are $1\le m<n\le4$, with leading coefficients $a,b$,
then $maDb-nbDa=0$ and $a^n/b^m=ch^k$. Homogeneity gives
$3k=4(n-m)$, so $(m,n)=(1,4)$.
Normalize $f=xh+a_4$ and $g=x^4+b_1x^3+b_2x^2+b_3x+b_4$.
The $x^3$ coefficient gives $D(a_4-hb_1/4)=0$; hence $a_4=hb_1/4$.
With $L=x+b_1/4$, equation \eqref{t:pairleading} becomes $Dg=0$, holding $L$
fixed. This source change preserves $R_{xx}=0$. The coefficient degrees
$1,2,3,4$ and \eqref{t:hinvariants} give
$f=Lh,g=L^4+\lambda Lh$.

For this exceptional nonbinary pair, subtract constant multiples of $R$
from $P,Q$ and choose univariate polynomials $\phi,\psi$ so that
$j=\max(\deg U,\deg V)$ is as small as possible in
\[
 P=xR+\phi(x)+U,\qquad Q=\psi(x)+V,
 \quad \deg\phi\le3,\quad \psi=x^4+\text{lower terms},\quad
 \deg U,\deg V\le3.
\]
If $j\ge1$, the degree-$(j+4)$ determinant equation is
\begin{equation}\label{t:residualpair}
 hD(V_j)-4x^3D(U_j)=0,
\end{equation}
where $U_j,V_j$ are the degree-$j$ parts. Indeed, subtracting $x$ times
the $R$ row from the $P$ row in the Jacobian leaves
$(R+\phi'+U_x,U_y,U_z)$; every omitted determinant term has lower degree.
But $\deg_xD(V_j)\le j-1\le2$, while the other term in
\eqref{t:residualpair} is divisible by $x^3$. Thus both vanish.
By \eqref{t:hinvariants}, $U_j,V_j\in\Span\{x^j\}$, with an additional
multiple of $h$ allowed only when $j=3$. Absorbing the $x^j$ terms into
$\phi,\psi$ and subtracting the corresponding multiples of $R$ lowers $j$,
a contradiction. Hence $U,V$ are constant and can also be absorbed:
$P=xR+\phi(x),Q=\psi(x)$, whose Jacobian with $R$ is zero.

Consequently $f,g$ are binary. If $a,b$ are the $x^3$ coefficients of
$P_3,Q_3$, the $x^2$ coefficient of $\Delta_7$ gives $D(ag-bf)=0$.
By \eqref{t:hinvariants} and independence of $f,g$, $a=b=0$. Write
\[
 P_3=x^2A+xA_2+A_3,\quad Q_3=x^2B+xB_2+B_3,
 \quad R_2=xC_1+C_2.
\]
There is no $x^2$ term in $R_2$ because $R_{xx}=0$.
The coefficients of $x$ in $\Delta_7$ and $x^3$ in $\Delta_6$ give,
with brackets in $(y,z)$,
\begin{align}
 [x]\Delta_7&=2\bigl(A\br gh-B\br fh\bigr)=0,\label{t:binary1}\\
 [x^3]\Delta_6&=-6h\br AB=0.\label{t:binary2}
\end{align}
Here the second identity follows from Euler's formula
$A\br Bh-B\br Ah=-3h\br AB$. If $(A,B)\ne(0,0)$, it gives
$(A,B)=\ell(\alpha,\beta)$ for a nonzero linear form $\ell$ and
a nonzero constant pair $(\alpha,\beta)$. Equation~\eqref{t:binary1} then
makes the nonzero quartic $\alpha g-\beta f$ a $D$-invariant,
contradicting \eqref{t:hinvariants}.
Therefore $A=B=0$. All terms of degree at least three are affine in $x$,
so $F_{xx}$ is constant and Proposition~\ref{t:directional} applies.
This completes the case $\rho=2$.

\subsection{Proof of tameness}\label{sec:tameness-assembly}
We first apply the leading-part calculations from
Section~\ref{sec:invertibility} to automorphisms. The cases with dependent
quartic leading components have been treated above.

\begin{proposition}\label{prop:tame-leading}
Let $F\in\Aut_{\C}(\A^3)$ have total degree at most four, and suppose
that the homogeneous quartic components of $F$ are linearly independent.
Then $F$ is tame.
\end{proposition}
\begin{proof}
The homogeneous factorization of Lemma~\ref{lem:factor} gives the four
patterns \eqref{eq:patterns}. In the two binary patterns, the proof of
Lemma~\ref{lem:binary-leading} produces a constant second derivative
in the $z$ direction. For the quadratic common-factor pattern, either the
same lemma applies or Lemma~\ref{lem:common-factor} contradicts the
Jacobian condition. The remaining pattern is a conic pencil: in the
absence of a square, Proposition~\ref{prop:square-forced} excludes it;
with a square, Proposition~\ref{prop:conic-complete} excludes the two
nonbinary normal forms and reduces the binary normal form to
Lemma~\ref{lem:binary-leading}. In every remaining case, the map therefore has a
constant second directional derivative. Proposition~\ref{t:directional}
then proves tameness.
\end{proof}

\begin{proof}[Proof of Theorem~\ref{thm:tame-aut}]
Let $\rho$ denote the dimension of the span of the homogeneous quartic
components. If $\rho=3$, apply Proposition~\ref{prop:tame-leading}.
If $\rho=0$, use the degree-three theorem of Blanc--van Santen
\cite[Corollary~3.9.5]{BVS}. For $\rho=1$, the two lower-degree components
are treated in Section~\ref{t:smallpair}; the only exceptional cubic-coordinate
case is excluded by Lemma~\ref{t:nagata}.
Finally, for $\rho=2$, choose a target component $R$ of degree at most
three. If it is of the form $x+C(y,z)$, Section~\ref{t:small} treats it
by local tameness and polynomial reduction. The remaining cubic normal
form is handled in Section~\ref{t:noncube}. These cases exhaust the
possibilities.
\end{proof}

\begin{proof}[Proof of Theorem~\ref{thm:combined}]
Theorem~\ref{thm:main} gives polynomial invertibility, and
Theorem~\ref{thm:tame-aut} gives tameness.
\end{proof}

\section{Concluding remarks}\label{sec:conclusion}
The two proofs address complementary aspects of polynomial maps in the
quartic degree range. The invertibility argument reduces the Keller
condition to homogeneous leading-part constraints and geometric
obstructions on cubic fibers. The tameness argument then combines those
constraints with coordinate classifications and local elementary
factorizations. The constant second directional derivative provides a
common criterion linking the two arguments.

The distinction between invertibility and tameness also affects the scope
of the field statements.
\begin{remark}[Scope of the field statement]
The descent of polynomial invertibility in Corollary~\ref{cor:charzero}
uses uniqueness of the formal inverse. It does not, by itself, descend a
tame factorization from $\C$ to an arbitrary characteristic-zero field.
Accordingly, Theorem~\ref{thm:combined} is formulated over $\C$.
\end{remark}

Two open questions provide natural directions for future study:
\begin{enumerate}[label=\textup{(\arabic*)}]
\item \textbf{Degree five.} Is every polynomial map
$F:\A^3_\C\to\A^3_\C$ of ordinary total degree five satisfying
$\det\J F\in\C^*$ a polynomial automorphism? If so, is it tame?
\item \textbf{Degree six.} Is every polynomial map
$F:\A^3_\C\to\A^3_\C$ of ordinary total degree six satisfying
$\det\J F\in\C^*$ a polynomial automorphism? If so, is it tame?
\end{enumerate}

Further directions include extracting explicit
tame factorizations from the reductions in Section~\ref{sec:tameness},
with quantitative bounds on the number and degrees of the elementary
factors. It is also natural to determine which parts of the leading-form and
residue arguments extend beyond the quartic degree range. Such extensions
would require a corresponding analysis of the additional leading-part
patterns and higher-degree coordinate components. The present proof also
suggests studying descent of tame factorizations separately from descent
of polynomial inverses.

\section*{Acknowledgements}
The authors used GPT-6 Astra for assistance with brainstorming,
mathematical development, and manuscript drafting. The authors are solely
responsible for the final content, analysis, and conclusions.

\section*{Funding}
Jie Wang was supported by the National Key R\&D Program of China under grant No.~2023YFA1009401 and the National Natural Science Foundation of China under grant No.~12171324.
Dan Yan was supported by the Scientific Research Fund of Hunan Provincial Education Department (Grant No. 25B0088), the National Natural Science Foundation of China (Grant No. 12371020), and the Construct Program of the Key Discipline in Hunan Province.

\end{document}